\RequirePackage{ifpdf}
\documentclass[11pt]{article}
\usepackage{amsthm}
\usepackage{amsmath}
\usepackage{amssymb}
\usepackage{subeqnarray}
\usepackage{cases}
\usepackage{latexsym}
\usepackage{graphicx}
\usepackage{fancybox,color}
\usepackage{enumerate,verbatim}
\usepackage{longtable}
\usepackage{color}
\usepackage{lscape}
\usepackage{graphicx}
\usepackage{float}
\definecolor{lred}{rgb}{1,0.8,0.8}
\definecolor{lblue}{rgb}{0.8,0.8,1}
\definecolor{dred}{rgb}{0.6,0,0}
\definecolor{dblue}{rgb}{0,0,0.5}
\usepackage{url}

\newcommand{\blue}[1]{\begin{color}{blue}#1\end{color}}

\newcommand{\red}[1]{\begin{color}{red}#1\end{color}}
\newcommand{\green}[1]{\begin{color}{green}#1\end{color}}

\def\cA{{\cal A}} 
\def\cB{{\cal B}}  \def\cF{{\cal F}}

\def\cD{{\cal D}}

\def\cI{{\cal I}}
\def\cP{{\cal P}}
 
\def\cS{{\cal S}}
\def\cT{{\cal T}} 
\def\cH{{\cal H}} 
\def\cG{{\cal G}}

\def\cS{{\cal S}} \def\S{{\cal S}}
\def\cY{{\cal Y}}
\def\cZ{{\cal Z}}

\newcommand{\inprod}[2]{\langle #1 , #2 \rangle}

\def\norm#1{\|#1\|}
\def\abs#1{\left| #1\right|}
\def\cC{{\cal C}}

\def\cQ{{\cal Q}}
\def\cK{{\cal K}}
\def\cI{{\cal I}}
\def\cW{{\cal W}}

\def\hL{\widehat{L}}

\def\inprod#1#2{\langle#1, \, #2\rangle}

\def\diag{{\rm diag}}

\def\cX{{\cal X}}

\newtheorem{theorem}{Theorem}[section]
\newtheorem{proposition}{Proposition}[section]
\newtheorem{assumption}{Assumption}[section]

\newtheorem{remark}{Remark}[section]

\newcommand{\A}{\mathcal{A}}

\newcommand{\D}{\mathcal{D}}

\newcommand{\T}{\mbox{\textrm{\tiny{T}}}}

\def\inprod#1#2{\langle #1,\,#2\rangle}
\def\norm#1{\|#1\|}

\DeclareMathOperator*{\argmini}{arg \, min}

\newcommand{\argmin}{\argmini}

\newcommand{\clear}[1]{}
\begin{document}

\title{\bf A Convergent 3-Block Semi-Proximal  Alternating Direction Method
of Multipliers for Conic Programming with $4$-Type of Constraints}

\author{
Defeng Sun\thanks{Department  of  Mathematics  and  Risk  Management  Institute, National University of Singapore, 10 Lower Kent Ridge Road, Singapore ({\tt matsundf@nus.edu.sg}). Research supported in part by the Ministry of Education, Singapore, Academic
Research Fund under Grant R-146-000-194-112.}, \,
Kim-Chuan Toh\thanks{Department of Mathematics, National University of Singapore, 10 Lower Kent Ridge Road, Singapore
({\tt mattohkc@nus.edu.sg}).
Research supported in part by the Ministry of Education, Singapore, Academic
Research Fund under Grant R-146-000-194-112.
}, \,
Liuqin Yang\thanks{Department of Mathematics, National University of Singapore, 10 Lower Kent Ridge Road, Singapore ({\tt yangliuqin@nus.edu.sg}).}
}
\date{April 9, 2014; Revised on October 21, 2014}
\maketitle

\begin{abstract}
The objective of this paper is  to design an efficient and convergent   alternating direction method of multipliers (ADMM)   for finding a solution of  medium accuracy to  conic programming problems whose constraints consist of linear equalities, linear inequalities,   a non-polyhedral cone and a polyhedral cone.   For this class of problems, one may apply the directly extended ADMM to their dual, which can be written in the form of  convex programming with four
separable blocks in the objective function and a coupling linear equation constraint. Indeed, the directly extended ADMM, though may diverge in theory,  often performs much better numerically  than many of its variants with theoretical convergence guarantee. Ideally,  one should find a convergent variant which is at least as efficient as the directly extended ADMM in practice. We achieve this  goal by  designing a  convergent semi-proximal ADMM (called sPADMM3c for convenience)     for   convex programming problems
having three separable blocks in the objective function with the third part  being linear.
 At each iteration, the proposed sPADMM3c takes  one special   block coordinate descent (BCD) cycle with the order $1 \rightarrow 3 \rightarrow 2 \rightarrow 3$, instead of the usual $1 \rightarrow 2 \rightarrow 3$ Gauss-Seidel BCD cycle used in the non-convergent  directly extended $3$-block ADMM,  for updating the  variable blocks. Our extensive numerical tests on the important class of
doubly non-negative  semidefinite programming (SDP) problems  with linear
equality and/or inequality constraints   demonstrate that our convergent method is at least {$20\%$} faster than the directly extended ADMM with unit step-length for
the vast majority of  about $550$ large scale
problems tested. This confirms  that at least for conic convex programming, one can design an ADMM with a special BCD cycle  of   updating the variable blocks
   can have both theoretical convergence guarantee and superior numerical efficiency over the directly extended
ADMM.

\end{abstract}



\section{Introduction}

Our primary motivation of this paper is to design  an efficient but   simple  first order method with guaranteed convergence to  find a solution of moderate accuracy   to the following conic programming with four types of constraints
\begin{eqnarray}
 ({\bf P}) \quad  \max \left\{ \inprod{-c}{x} \mid \cA_{{E}} x=  b_{ E}, \; \cA_I x \ge b_I ,\; x \in \cK, \;   x\in \cK_p \right\}.
  \nonumber 
\end{eqnarray}
Here $\cA_E$ and $\cA_I$ are two linear maps defined from the finite-dimensional real Euclidean space
 $\cX$ to $\Re^{m_E}$ and $\Re^{m_I}$, respectively, $(b_E, b_I) \in \Re^{m_E}\times  \Re^{m_I}$  and $c\in \cX$ are given data,
$\cK\subseteq \cX$ is a pointed closed convex (non-polyhedral) cone whose interior ${\rm int}(\cK) \neq \emptyset$ and $\cK_p$ is a polyhedral convex cone in $\cX$
such that $\cK \cap \cK_p$ is non-empty.
Note that in theory the fourth block constraint in ({\bf P})  can be incorporated into the first and the second block constraints. However,
 treating the polyhedral cone $\cK_p$ separately instead of representing it in terms of  equalities and inequalities is of considerable advantage in numerical computations.
  Problem ({\bf P}) covers a wide range of  interesting convex models.  In particular, it includes the important class of doubly non-negative (DNN) semidefinite programming (SDP) with both equality and inequality constraints
\begin{eqnarray}
{\bf (SDP)} \quad  \max \left\{ \inprod{-C}{X} \mid\cA_E X=  b_{ E}, \; \cA_I X \ge b_I ,\; X \in \cS_+^n, \; X\in \cK_p \right\},
 \nonumber 
\end{eqnarray}
where $\cS_+^n$ is the cone of $n\times n$ symmetric and positive semidefinite matrices in the space  of   $n\times n$ symmetric matrices $\cS^n$, $X\in \cK_p$ means that
every entry of the matrix  $X\in \cS^n$ is non-negative (one can, of course, only require
a subset of the entries of $X$ to be non-negative or non-positive or free) and $C\in \cS^n$ is a given symmetric matrix.

For a given linear map $\cA$ from $\cX$ to another   finite-dimensional real Euclidean space  $\cX'$, we denote its adjoint by $\cA^*$ and
 for any closed convex set $\cC\subseteq \cX$, we denote the metric projection operator onto $\cC$ by $\Pi_\cC(\cdot)$. If   $\cC \subseteq \cX$ is a   closed convex cone, we denote its dual cone   by $\cC^*$, i.e.,
 \[
 \cC^* : =\{ d \in \cX \, |\, \langle d, x\rangle \ge 0 \quad  \forall\, x\in \cC\}.
 \]
 In this paper, we will make extensive use of the Moreau decomposition theorem \cite{moreau1962decomposition}, which states that $x = \Pi_\cC(x) - \Pi_{\cC^*}(-x)$ for any nonempty closed convex cone $\cC \subseteq \cX$ and    $x\in \cX$.

   The dual of ({\bf P}) takes the form of
\begin{eqnarray}
  \min \left\{ -\inprod{b_I}{y_I} - \inprod{b_E}{y_E}\mid s+\cA_I^* y_I +z+ \cA^*_E y_E = c,\;
   s\in \cK^*, \; z\in\cK^*_p ,\; y_I \ge 0\right\},
  \label{eq-dualconic0}
\end{eqnarray}
which can equivalently be written as the following convex programming with four separate blocks in the objective function and a coupling linear equation constraint:
\begin{eqnarray}
 ({\bf D})\quad  \min \left\{ \delta _{{\cK}^*} (s) +( \delta _{\Re_+^{m_I}} (y_I)    -\inprod{b_I}{y_I})+\delta_{\cK^*_p}(z) - \inprod{b_E}{y_E}   \mid s+\cA_I^* y_I +z+ \cA^*_E y_E = c
     \right\},
  \nonumber 
\end{eqnarray}
where for any given  set $\cC$, $\delta_\cC(\cdot)$ is the indicator function over $\cC$ such that $\delta_\cC (u)=0$ if $u\in \cC$ and $\infty$ otherwise.

Problem  ({\bf D}) belongs to a class of multi-block convex optimization problems whose objective function is the sum of  $q$ convex functions without overlapping variables:
\begin{eqnarray}
 \min \left\{ {{\sum}}_{i=1}^q\theta_i(w_i)  \mid {\sum}_{i=1}^q\cB^*_i w_i = c \right\},
  \label{eq-ADMMq}
\end{eqnarray}
where for each $i\in \{1, \ldots, q\}$, $ \cW_i$ is a  finite dimensional real Euclidean space equipped with an inner product $\inprod{\cdot}{\cdot}$ and its induced norm $\norm{\cdot}$,  $\theta_i: \cW_i \mapsto (-\infty,+\infty]$ is a closed proper convex functions, $\cB_i:   \cX \mapsto \cW_i$ is a  linear map and  $c \in \cX$ is given. Note that one  can write ({\bf D}) in the form  of (\ref{eq-ADMMq}) in a number of different ways. One natural choice is of course to write  ({\bf D}) in terms of (\ref{eq-ADMMq}) for $q=4$
 with $(w_1, w_2, w_3, w_4) \equiv (s, y_I, z, y_E)$. However, by noting that in   ({\bf D}) the objective function containing the $y_E$ part is a linear term, we shall treat  ({\bf D})  as a special case of  (\ref{eq-ADMMq}) for $q=3$ with
  $(w_1, w_2, w_3) \equiv (s, y_I, (z, y_E))$. In the latter case, the third function $\theta_3 $ is partially linear: it is linear about $y_E$ but nonlinear about $z$.
  This partial linear structure of $\theta_3$ will be
heavily
exploited in our pursuit of designing a convergent but efficient  first order method  in this paper.

Let $\sigma > 0$ be given. The augmented Lagrange function for (\ref{eq-ADMMq}) is defined by
\begin{eqnarray}
L_\sigma(w_1,\ldots,w_q;x) :=  {\sum}_{i=1}^q\theta_i(w_i) + \inprod{x}{ {\sum}_{i=1}^q\cB^*_i w_i - c}+ \frac{\sigma}{2} \norm{ {\sum}_{i=1}^q\cB^*_i w_i  - c}^2\nonumber
\label{eq-ADM-q-AL}
\end{eqnarray}
for $w_i\in \cW_i$, $i=1, \ldots, q$ and $x\in \cX.$
%
Choose any initial points $w_i^0 \in {\rm dom}(\theta_i) $, $i=1, \ldots, q$ and $x^0 \in \cX$.  The classic augmented Lagrangian function method of Hestenes-Powell-Rockafellar
\cite{Hestenes1969,Powell1969,rockafellar1976augmented}
consists of the following iterations:
\begin{eqnarray}
 (w_1^{k+1}, \ldots, w_q^{k+1}) &= & \argmin
L_\sigma(w_1,\ldots,w_q;x^k),
\label{eq-auglagl-w}
\\[5pt]
 x^{k+1} &=&\displaystyle  x^k +\tau \sigma\left ({\sum}_{i=1}^q\cB^*_i w_i^{k+1}  -c\right),
\label{eq-auglagl-x}
\end{eqnarray}
where $\tau >0$, e.g., $ \tau\in(0,2)$,  is a positive constant that controls the step-length in \eqref{eq-auglagl-x}. To solve (\ref{eq-auglagl-w}) exactly  or approximately
to high precision can be a challenging task in many situations.  To deal with this challenge, one may  try to replace \eqref{eq-auglagl-w}  by considering the following
$q$-block alternating direction method of multipliers (ADMM):
\begin{eqnarray}\label{admmlb}
 w_1^{k+1} &=& \argmin
  L_\sigma(w_1,w_2^k\ldots,w_q^k;x^k), \nonumber
\\[-8pt]
& \vdots & \nonumber
\\[0pt]
 w_i^{k+1} &=& \argmin L_\sigma(w_1^{k+1},\ldots, w_{i-1}^{k+1}, w_i, w_{i+1}^k, \ldots, w^k_q;x^k), \nonumber
 \\[-8pt]
& \vdots &
\\[0pt] \nonumber
  w_q^{k+1} &=& \argmin L_\sigma(w_1^{k+1},\ldots,   w_{q-1}^{k+1},  w_q;x^k), \nonumber
\\[5pt]
 x^{k+1} &=& \displaystyle  x^k +\tau \sigma\big ({{\sum}_{i=1}^q}\cB^*_i w_i^{k+1}  -c\big).
 \nonumber
\end{eqnarray}
The above $q$-block ADMM   is  extended directly from the ADMM
 for solving the following
  2-block convex optimization problem
\begin{eqnarray}
 \min \left\{ \theta_1(w_1)+\theta_2(w_2) \mid \cB_1^* w_1 + \cB^*_2 w_2   = c \right\}.
  \label{eq-ADMM2}
\end{eqnarray}
For a chosen initial point $(w_1^0, w_2^0,  x^0) \in {\rm dom}(\theta_1)\times  {\rm dom}(\theta_2) \times \cX$,  the classic $2$-block ADMM  consists of the iterations:
\begin{eqnarray}
 w_1^{k+1} &=& \argmin
   L_\sigma(w_1,w_2^k; x^k),
 \label{eq-ADM-x}
\\[5pt]
 w_2^{k+1} &=& \argmin L_\sigma(w_1^{k+1},w_2;x^k),
 \label{eq-ADM-y}
\\[5pt]
x^{k+1} &=& x^k +  \tau \sigma (\cB^*_1 w_1^{k+1} + \cB_2^* w_2^{k+1} -c),
 \label{eq-ADM-w}
\end{eqnarray}
where $\tau  >0$ is a positive constant.
The classic  $2$-block ADMM for solving \eqref{eq-ADMM2} was first introduced by  Glowinski and
Marrocco \cite{glowinski1975sur} and Gabay and Mercier \cite{gabay1976dual}.
When ${\cal B}_1^*=\cI$, the identity mapping, $\cB^*_2$ is injective and $\theta_1$ is strongly convex, the
convergence of the classic $2$-block ADMM has been proven first by Gabay and Mercier \cite{gabay1976dual}
 for any $\tau \in(0,2 )$ if  $\theta_2$ is linear, and then
by Glowinski \cite{Glowinski1980lectures} and Fortin and Glowinski \cite{fortin1983augmentedlag}   for any  $\tau \in (0, (1+\sqrt{5}\,)/2)$ if $\theta_2$ is a general nonlinear convex  function. Gabay  \cite{Gabay1983299} has further shown  that
   the classic $2$-block ADMM for $\tau =1$ is a special case of the Douglas-Rachford splitting method. Moreover, Eckstein and Bertsekas \cite{eckstein1992douglas} has  shown that the latter is actually an application of the proximal point algorithm on the dual problem by means of a specially-constructed splitting operator.
A variation of the classic $2$-block ADMM is to adjust the penalty parameter  $\sigma$ at every iteration based on the previous iterations'  progress with the goal of improving the convergence in practice. That is, one replaces (\ref{eq-ADM-w}) by
\begin{eqnarray}
x^{k+1} = x^k + \tau \sigma_k (\cB^*_1 w_1^{k+1}+\cB^*_2 w_2^{k+1}-c).
\end{eqnarray}
A scheme to adjust the penalty parameters $\sigma_k$  was studied in \cite{he2000alternating,wang2001decomposition}, which often works well in practice.
 Due to its extreme simplicity and efficiency in several applications in mathematical imaging science, signal processing and etc., the classic  $2$-block ADMM has regained its popularity in recent years.
 For a tutorial on the classic $2$-block ADMM, one may consult the  the recent work by Eckstein and Yao
\cite{eckstein2012}.

The multiple-block ADMM  with larger $\tau$  often
works very well in many cases. For example, Wen, Goldfarb and Yin \cite{wen2010alternating}   used  the $3$-block  ADMM with $\tau =1.618$ to design an efficient  software for  solving some SDP problems  of large sizes.  However, it is shown very recently  by Chen, He, Ye, and Yuan \cite{chen2013direct} that  in contrast to the classic 2-block ADMM, the   directly  extended 3-block ADMM  with $\tau =1$
may diverge\footnote{The final version of \cite{chen2013direct} includes  a non-convergent  example for the ADMM with any prefixed $\tau$ not smaller than  $10^{-8}$ and the penalty parameter $\sigma= 1$.}.  This dashes any hope of using the directly extended $q$-block ($q\ge 3$) ADMM   without
modifications\footnote{Hong and Luo \cite{HongL2012} provided a proof on the convergence of the directly extended  $q$-block ADMM under some restrictive assumptions including a global error bound condition with a sufficiently small step-length $\tau$. Since in practical computations, one always prefers a larger step-length for better numerical efficiency, a convergence result of this nature is mainly of theoretical importance.}.
 Actually, even before the announcement of \cite{chen2013direct},
several researchers have made serious  attempts in correcting the possible non-convergence of the multiple-block ADMM \cite{HYZCai2013, HTXYuan2013, he2012alternating,  HTYuan2014, HYuan2013,TYuan2011}, to name only a few. A recent work by Wang, Hong, Ma and Luo \cite{WHMLuo2013} is also along this line.  Among the work on correcting the non-convergence, the $q$-block ADMM with an additional  Gaussian back substitution \cite{he2012alternating} distinguishes itself for its simplicity and generality. However, to the best of our knowledge, up to now
 the dilemma is that at least for convex conic programming,  the modified versions though with convergence guarantee,
often perform $2$--$3$ times slower than the multi-block ADMM with no convergent guarantee.

In this paper we aim to resolve the dilemma just mentioned
by focusing on the conic programming problem ({\bf P}).
We achieve this  goal by  proposing a  convergent  semi-proximal ADMM  (sPADMM)
first  for   convex programming problems having three separate blocks in the objective function with the third part  being linear (we call this method sPADMM3c) and then extend it to the general case.  Our extensive numerical tests on the important class of
doubly non-negative  semidefinite programming (SDP) problems  with linear
equality and/or inequality constraints   demonstrate that our convergent method is at least {$20\%$} faster than the directly extended ADMM with unit step-length for
the vast majority of  about $550$ large scale
problems tested. This confirms  that our convergent
sPADMM
can have both theoretical convergence guarantee and superior numerical efficiency over the directly extended ADMM.

Note that \eqref{admmlb} performs  the  usual $1 \rightarrow 2 \rightarrow \cdots \rightarrow q$ Gauss-Seidel block coordinate descent (BCD) cycle 
in an
alternating
way of minimizing the   variable blocks
  for solving \eqref{eq-auglagl-w} inexactly.
   In contrast, the proposed sPADMM3c
which will be presented in Section \ref{subsection:h-linear} for the case $q=3$,  takes  the  special   $1 \rightarrow 3 \rightarrow 2 \rightarrow 3$ BCD cycle at each iteration.
This special cycle  actually uses an  essentially BCD cyclic rule in the terminology of Tseng \cite{tseng2001convergence} to minimize the   variable blocks for solving \eqref{eq-auglagl-w} inexactly.  Given the fact that all the  three component functions in the objective  of a counterexample constructed in \cite{chen2013direct} to illustrate the non-convergence of the directly extended 3-block ADMM are zero functions, it comes as a pleasant surprise that one only needs to update  the third variable block twice to get a convergent $3$-block ADMM provided that $\theta_3$ is linear. At the moment,  it  is not clear to us if our  corresponding ADMM  is still  convergent when $\theta_3$ is nonlinear though we conjecture that it is true\footnote{In a recent paper \cite{li2014schur},  Li, Sun and Toh  report that the conjecture is indeed true if   $\theta_3$ is a convex quadratic function.}. In any case, for the conic programming problem   ({\bf P}),  the requirement on the linearity of $\theta_3$ is not
restrictive at all, as 
we will see in the subsequent analysis.

The remaining parts of this paper are organized as follows. In the next section, for our subsequent developments we will present in details the convergent properties of a semi-proximal ADMM for solving two-block convex optimization problems. In Section \ref{section:PADMM3C}, we will   introduce a convergent
$3$-block
sPADMM
first  for the special case where  the third function $\theta_3$ is a linear function and then show how this approach can be applied to the case where $\theta_3$ is partially linear or fully nonlinear.  We should emphasize  that  the linear or partial linear structure of $\theta_3$ will not render the directly extended $3$-block ADMM to become convergent as  the three functions $\theta_1,\theta_2$ and $\theta_3$ constructed  in the  counterexample  in   \cite{chen2013direct} are all linear (actually, zero functions).
Extensions to the multi-block case are also briefly discussed in this section.
In Section \ref{section:ConicProgramming}, the applications of the convergent $3$-block
sPADMM
to conic programming
 are discussed.  Section \ref{section:numericalresults} is
devoted to the implementation
and numerical experiments of using our convergent $3$-block
sPADMM
for solving ({\bf SDP}). We separate our numerical tests into two parts.
The first part is on the  doubly non-negative  SDP, i.e., ({\bf SDP}) without the inequality constraints ${\cal A}_I X\ge b_I$. For this class of  SDP problems,   our proposed convergent $3$-block ADMM is
more  efficient than the directly extended $3$-block ADMM (with $\tau =1$) and it is competitive to the latter with $\tau=1.618$
 in terms of both the number of iterations and computing time.
We should mention again that the directly extended $3$-block ADMM
has no convergent guarantee regardless of whether $\tau=1$ or $\tau =1.618$.
The second part of the section is on a class of general    ({\bf SDP}) having four types of constraints  including a  large number of inequality constraints  ${\cal A}_I X\ge b_I$. This time,
our proposed convergent $3$-block
sPADMM
 is still more  efficient than the directly extended
$4$-block
sPADMM with $\tau =1$ and is competitive to the latter even for $\tau=1.618$  in terms of the  computing time and needs less number of iterations.
Our research conducted here opens up the possibility of designing an efficient and convergent  ADMM with a suitably chosen essentially BCD cyclic rule
rather than  the usual  Gauss-Seidel BCD cyclic  rule where the latter may lead to a non-convergent ADMM,
for multi-block convex optimization problems
(\ref{eq-ADMMq}) with structures beyond those considered in ({\bf SDP}). We conclude our paper in the final section.

\bigskip
\noindent
{\bf Notation.}
For any given self-adjoint positive semidefinite operator $\cT$ that maps a real Euclidean space $\cX$ into itself, we let  $ \cT^{1/2}$ be the unique   self-adjoint positive semidefinite operator
such that $\cT^{1/2} \cT^{1/2} =\cT $ and define
\[
\|x\|_{\cT} : = \sqrt{\inprod{x}{\cT x}} =\|\cT^{1/2} x\|\quad \forall\, x\in \cX.
\]
For a given convex function $\phi: \cX \to (-\infty, \infty]$, $\phi^*$ denotes its Fenchel conjugate, i.e.,
\[
\phi^* (x) : =\sup\{\inprod{y}{x} -\phi(y)\}, \quad x\in \cX.
\]


\section{Preliminaries}
\label{section:ADMM2B}
Thoughout this paper, we assume that  $ \cX,  \cY, \cZ$ are three finite dimensional real Euclidean spaces each equipped with an inner product $\inprod{\cdot}{\cdot}$ and its induced norm $\norm{\cdot}$. Let
$f: \cY \mapsto (-\infty,+\infty]$,  $g: \cZ \mapsto (-\infty,+\infty]$   be given closed proper convex functions, $\cF: \cX \mapsto \cY$, $\cG: \cX \mapsto \cZ$ be given linear maps and  $c \in \cX$ be also given.
Let $\cS$ and $\cT$ be two given self-adjoint positive semidefinite
(not necessarily positive definite) linear operators on $\cY$ and $\cZ$,
respectively.

The purpose of this section is to discuss the convergent properties of a
semi-proximal ADMM (sPADMM)
for solving the 2-block convex optimization problem
\begin{eqnarray}
\min \left\{ f(y)+g(z) \mid \cF^* y + \cG^* z = c \right\}
\label{eq-ADMM2-nonlinear}
\end{eqnarray}
and its dual
\begin{eqnarray}
\max \left\{ -\inprod{c}{x}-f^*(-\cF x)-g^*( -\cG x )  \right\},
\label{eq-ADMM2-primal0}
\end{eqnarray}
which can equivalently  be written as
\begin{eqnarray}
\min \left\{\inprod{c}{x} + f^*(u)+g^*(v)\mid  \cF x +u =0, \; \cG x + v=0 \right\}.
\label{eq-ADMM2-nonlinear-primal}
\end{eqnarray}
These convergence results are the essential ingredients in proving the convergence of our $3$-block
sPADMM proposed in Section \ref{section:PADMM3C}.

Recall that the augmented Lagrange function  for problem  \eqref{eq-ADMM2-nonlinear} is defined by
\begin{eqnarray}
L_\sigma(y,z;x) = f(y) + g(z) + \inprod{x}{\cF^* y + \cG^* z  - c} + \frac{\sigma}{2} \norm{\cF^* y + \cG^* z - c}^2.
\label{eq-ADM-AL}
\end{eqnarray}
It is clear that in order for the classic 2-block ADMM scheme \eqref{eq-ADM-x}--\eqref{eq-ADM-w}    applied to problem  \eqref{eq-ADMM2-nonlinear} to work, one needs to assume that both  subproblems have a solution. The existence of solutions for the subproblems can be guaranteed if we assume that  the objective functions in (\ref{eq-ADM-x}) and (\ref{eq-ADM-y}) are both coercive. However,  conditions ensuring the boundedness of the two generated sequences $\{y^{k+1}\}$ and $\{z^{k+1}\}$ are very subtle while the boundedness of the dual variable sequence $\{x^{k+1}\}$ is readily obtainable. More importantly, it is also desirable that  both  $y^{k+1}$ and $z^{k+1}$ can be computed relative easily  if $f$ and $g$ have conducive structures.  In this regard, the  following
sPADMM
is preferred.

\bigskip
\centerline{\fbox{\parbox{\textwidth}{
			{\bf Algorithm sPADMM2}: {\bf A generic 2-block semi-proximal ADMM for solving \eqref{eq-ADMM2-nonlinear}.}
			\\[5pt]
			Let   $\sigma > 0$ and  $\tau \in (0,\infty)$ be given parameters.
			Choose $y^0 \in {\rm dom}(f)$, $z^0 \in {\rm dom}(g)$, and
			$x^0 \in \cX$. Perform the $k$th iteration as follows:
			\begin{description}
				\item[Step 1.]  Compute
				$\displaystyle y^{k+1} = \argmin L_{\sigma}(y,z^k; x^k)+\frac{\sigma}{2}\norm{  y-y^k}^2 _{\cS}.$
				\item[Step 2.]  Compute
				$\displaystyle z^{k+1} = \argmin  L_\sigma(y^{k+1},z; x^k)+\frac{\sigma}{2}\norm{ z-z^k}^2 _{\cT}.$
				\item[Step 3.] Compute $\displaystyle x^{k+1} = x^k + \tau \sigma (\cF^*y^{k+1} + \cG^*z^{k+1}-c)$.
			\end{description}
		}}}
		\bigskip
		
		In the above 2-block
sPADMM for solving problem \eqref{eq-ADMM2-nonlinear}, the choices of $\cS$ and $\cT$ are very much problem dependent. The general principle is that both $\cS$ and $\cT$ should be as small as possible while $y^{k+1}$ and $z^{k+1}$ are still relatively easy to compute.  The convergence analysis of the 2-block
sPADMM can be conducted by following the proof given by  Fortin and Glowinski \cite{fortin1983augmentedlag} based on variational analysis. This has been done in
		\cite{fazel2013hankel}. Here we will only summarize what we need for later developments. For details, one may refer to Appendix B in \cite{fazel2013hankel}.
	
		\bigskip
		\begin{assumption}\label{assumption:CQ2}
			There exists $(\hat y,\hat z)\in {\rm ri}({\rm dom}\,f\times{\rm dom}\,g) $ such that $\cF^* \hat y + \cG^* \hat z=c$.
		\end{assumption}

		Under Assumption \ref{assumption:CQ2}, it follows from \cite[Corollary~28.2.2]{rockafellar1970convex} and \cite[Corollary~28.3.1]{rockafellar1970convex} that $(\bar{y},\bar{z}) \in  \cY \times \cZ$ is an optimal solution to problem \eqref{eq-ADMM2-nonlinear} if and only if there exists a Lagrange multiplier $\bar{x} \in \cX$ such that
		\begin{equation}\label{genopt2B}
		0\in \cF \bar{x} + \partial f(\bar{y}), \quad 0 \in \cG \bar{x} + \partial
		g(\bar{z}),\quad   \cF^*\bar{y}+\cG^*\bar{z} -c=0,
		\end{equation}
		where $\partial f$ and $\partial g$ are the subdifferential mappings of $f$ and $g$ respectively. Moreover, any $\bar{x}\in \cX$ satisfying \eqref{genopt2B} is an optimal solution to the dual problem \eqref{eq-ADMM2-nonlinear-primal}.		
		Since both  $\partial f$ and $\partial g$  are maximal monotone \cite[Theorem~12.17]{rockafellar1998variational}, there exist two
		self-adjoint and positive semidefinite operators  $\Sigma_f$ and $\Sigma_g $    such that   for all $y,  y'\in {\rm dom}(f)$, $u\in \partial f(y )$,
		and ${u}'\in \partial f({y}')$,
		\begin{equation}\label{monosub1}
		\langle u-{u}', y -{y}'\rangle \geq
		\|y -{y}'\|^2_{\Sigma_f}
		\end{equation}
		and for all $z,  z' \in {\rm {dom}}(g)$, $v\in \partial
		g(z )$, and $ {v}'\in \partial g({z}')$,
		\begin{equation}\label{monosub2}
		\langle v-{v}', z -{z}'\rangle \geq
		\|z -{z}'\|^2_{\Sigma_g}.
		\end{equation}

		\begin{theorem}\label{thempadm}  Let
			$\Sigma_f$ and $\Sigma_g $    be   the   two self-adjoint and positive semidefinite operators defined by (\ref{monosub1}) and (\ref{monosub2}), respectively.
			Suppose  that the solution set of problem \eqref{eq-ADMM2-nonlinear} is nonempty and that  Assumption \ref{assumption:CQ2}  holds.   Assume that
			$\cS$ and $\cT$ are chosen such that the
			sequence $\{(y^k,z^k,x^k)\}$  generated by Algorithm sPADMM2 is well defined. Let $(\bar y, \bar z)$ be any  optimal solution to problem \eqref{eq-ADMM2-nonlinear} and $\bar{x}$ be any optimal solution to problem \eqref{eq-ADMM2-nonlinear-primal}, respectively.
			For $k=0,1, 2, \ldots, $ denote
				\[
			y_e^k: =y^k -\bar y, \quad z_e^k: =z^k -\bar z \quad {\rm and}\quad x_e^k: = x ^k -\bar {x}.
			\]
			Then, under the condition either
			(a)  $\tau\in (0, (1+\sqrt{5}\,)/2)$ or (b)  $\tau \ge (1+\sqrt{5}\, )/2$ but $\sum_{k=0}^\infty (\|\cG^*(z^{k+1}-z^k)\|^2 + \tau^{-1} \|\cF^* y^{k+1} +\cG^* z^{k+1}-c\|^2) <\infty$,
			the following results hold:
			\begin{enumerate}
				\item [{\rm (i)}] The sequence $\{\|x^{k+1}_e\|^2+\|z^{k+1}_e\|^2 _{(\sigma^{-1} \Sigma_g + \cT + \cG \cG^*)}+\|y^{k+1}_e\|^2 _{(\sigma^{-1}\Sigma_f +\cS +  \cF\cF^*)}\}$  is bounded.
				\item [{\rm (ii)}]   If $(y^\infty, z^\infty, x^\infty)$ is an accumulation point of $\{(y^k, z^k, x^k)\}$, then $(y^\infty, z^\infty)$ solves \eqref{eq-ADMM2-nonlinear} and $x^\infty$  solves \eqref{eq-ADMM2-nonlinear-primal}, respectively, and it holds that
				\[
				\lim_{k\rightarrow \infty}\left(  \|x^{k+1}_e\|^2+\|z^{k+1}_e\|^2 _{(\sigma^{-1}\Sigma_g +\cT + \cG\cG^*)}+\|y^{k+1}_e\|^2 _{(\sigma^{-1}\Sigma_f +\cS +  \cF \cF^*)}\right)=0,
				\]
				where in the definition of  $(y_e^k, z^k_e, x_e^k)$, the point $(\bar y, \bar z, \bar x)$ is replaced by $(y^\infty, z^\infty, x^\infty)$.
				
				\item [{\rm (iii)}] If both $\sigma^{-1}\Sigma_f +\cS +  \cF\cF^*$ and $\sigma^{-1}\Sigma_g + \cT +\cG \cG^*$ are positive definite, then  the  sequence
				$ \{(y^k, z^k, x^k)\}$, which  is automatically well defined, converges to a unique limit, say, $(y^\infty, z^\infty, x^\infty)$  with $(y^\infty, z^\infty)$ solving \eqref{eq-ADMM2-nonlinear} and $x^\infty$ solving  \eqref{eq-ADMM2-nonlinear-primal}, respectively.
				\item [{\rm (iv)}]  When the $z$-part disappears, the corresponding results in parts (i)--(iii) hold under the condition either $\tau \in (0,2)$ or $\tau \ge2$ but $\sum_{k=0}^\infty \|\cF^* y^{k+1} -c\|^2<\infty$ .
			\end{enumerate}
		\end{theorem}

		\begin{remark}\label{rem:convergencePADMM2}
			The conclusions of Theorem \ref{thempadm} for the case that $\tau \in (0, (1+\sqrt{5}\,)/2)$ follow directly from the results given in  \cite[Theorem B.1]{fazel2013hankel}.  For the case that  $\tau \ge (1+\sqrt{5}\,)/2$ but $\sum_{k=0}^\infty (\|\cG^*(z^{k+1}-z^k)\|^2 + \tau^{-1} \|\cF^* y^{k+1} +\cG^* z^{k+1}-c\|^2) <\infty$, we can just mimick the proofs for  {part (c)} in  \cite[Theorem B.1]{fazel2013hankel} for the case that $\tau \in (1, (1+\sqrt{5}\,)/2)$ by using
			{part (b)} in  \cite[Theorem B.1]{fazel2013hankel} and the property on  Fej\'er monotone sequences. Similarly, the conclusions for part (iv) can be derived correspondingly by using part (d) in \cite[Theorem B.1]{fazel2013hankel}. In our numerical computations, we always start with  a larger $\tau$, e.g., $\tau =1.95$, and reset it as $\tau :=\max(\rho \tau, 1.618)$ for some $\rho\in (0,1)$ if at the $k$-th iteration
			\[
			\|\cG^*(z^{k+1}-z^k)\|^2 + \tau^{-1} \|\cF^* y^{k+1} +\cG^* z^{k+1}-c\|^2> c_0k^{-1.2}
			\]
			for some constant $c_0>0$. Since $\tau$ can  be reset for a finite number of times only,  eventually either condition (a) or condition (b) in Theorem \ref{thempadm} is satisfied. Consequently, the conclusions of parts (i)-(iii) in  Theorem \ref{thempadm} hold.
			When the    $z$-part disappears, we can start with $\tau \geq 2  $  and reset it accordingly by using a similar procedure to the above.
		\end{remark}

		\begin{remark}
		Independent of Fazel et al. \cite{fazel2013hankel}, Deng and Yin \cite{dengglobal} also analyze the global convergence of Algorithm sPADMM2, though the focus of \cite{dengglobal} is mainly on analyzing the  rate of convergence, for the following   cases: i) $\cS \succeq 0$, $\cT \succeq 0$,  $\tau \in (0,1]$;  ii) $\cS \equiv 0$, $\cT \succeq 0$,   $\tau \in [1,\frac{1+\sqrt{5}}{2})$. The most interesting case used in  this paper  of taking  $\cS \succneqq  0$, $\cT \succneqq 0$  and   $\tau >1$ (in particular, $\tau =1.618$) is not covered by \cite{dengglobal}.
		\end{remark}

\section{A convergent $3$-block semi-proximal ADMM}\label{section:PADMM3C}
  Assume that $\cW$ is a finite dimensional real Euclidean space.  Let
 $h: \cW \mapsto (-\infty,+\infty]$  be a given closed proper convex function,  $\cH: \cX \mapsto \cW$ be a given   linear map. For the subsequent discussions, we let
 $\cT_f$ and $\cT_g$ be two given self-adjoint positive semidefinite
 (not necessarily positive definite) linear operators on $\cY$ and $\cZ$,
 respectively.

 Consider the following $3$-block convex optimization problem
  \begin{eqnarray}
 \min \left\{ f(y)+g(z)+h({w}) \mid \cF^* y + \cG^* z + \cH^* w = c \right\}.
 \label{eq-ADMM3-linear}
\end{eqnarray}
The dual of \eqref{eq-ADMM3-linear} is
\begin{eqnarray}
   \max \left\{ -\inprod{c}{x}-f^*(-\cF x)-g^*( -\cG x )-h^*(-\cH x)  \right\},
  \label{eq-ADMM3-primal0}
\end{eqnarray}
which can equivalently  be written as
\begin{eqnarray}
   \min \left\{\inprod{c}{x} + f^*(u)+g^*(v) +h^*(s)\mid \cF x +u =0, \; \cG x + v=0,\;  \cH x +s = 0 \right\}.
  \label{eq-ADMM3-linear-primal}
\end{eqnarray}
By noting that the three variables $u$, $v$, and $s$ are decoupled in the constraints of problem \eqref{eq-ADMM3-linear-primal},  one may attempt to apply the classic $2$-block ADMM if
 $\cF^* \cF+\cG^* \cG +\cH^*\cH$ is positive definite or the 2-block
sPADMM if it  is only positive semidefinite, to (\ref{eq-ADMM3-linear-primal}) with $x$ and $(u,v,s)$ as two separate block variables.\footnote{This comment can be direcly applied to the $q$-block convex optimization problem (\ref{eq-ADMMq}).}
   However, as far as we know from our numerical experiments,
this approach is    less efficient than working with the problem of the form \eqref{eq-ADMM3-linear} directly.
    A possible explanation for this phenomenon is that
   the sizes of  $ \cF^* \cF+\cG^* \cG +\cH^*\cH $ are often too large to admit an efficient Cholesky factorization  and consequently one is forced  to add a large semi-proximal term to it to  make
the new operator more amenable for
practical  computations. At least, this is the case for conic programming ({\bf P}).


In the next subsection, we will introduce our approach first  for the   case where $h$ is a linear function and prove the convergence of our approach by relating it to a particularly designed 2-block
sPADMM.
For ({\bf D}), this  corresponds to the case where $m_I= 0$ or $\cK_p=\cX$. After that, we will extend our idea to deal with  the case where $h$ is only partially linear or fully nonlinear and the $q$-block case in Section \ref{section:extension}.

\subsection{The case where  $h$ is linear}\label{subsection:h-linear}
In this subsection, we are particularly interested in the case where $h$ is a linear function of the form
  \begin{eqnarray}
h(w): =-\inprod{b}{w} \quad \forall\, w\in \cW,
\label{eq:linearmappingh}
\end{eqnarray}
where $b\in \cW$ is given.  For simplicity, by removing the redundancy if necessary (although
it may not be an easy task numerically), we  assume that $\cH\cH^*$ is invertible, i.e., $\cH$ is surjective.

 For a given $\sigma>0$, let $L_\sigma(y,z,w;x)$ be the augmented Lagrange function for \eqref{eq-ADMM3-linear}, i.e., for any $(y,z,w,x)\in \cY\times \cZ\times \cW\times \cX$,
 \begin{eqnarray}
 L_\sigma(y,z,w;x) =  f(y)+g(z)-\inprod{b}{w} + \inprod{x}{ \cF^* y + \cG^* z + \cH^* w - c}+ \frac{\sigma}{2}
 \norm{ \cF^* y + \cG^* z + \cH^* w  - c}^2.
 \label{eq-ADM-3-AL}
 \end{eqnarray}
The following constraint qualification is needed for our subsequent discussions.

\begin{assumption}\label{assumption:CQ1}
 There exists $(\hat y,\hat z, \hat w)\in {\rm ri}({\rm dom}\,f\times{\rm dom}\,g )\times \cW $ such that $\cF^* \hat y + \cG^* \hat z +\cH^* \hat{w}=c$.
\end{assumption}

Similar to the discussion in Section \ref{section:ADMM2B}, under Assumption \ref{assumption:CQ1}, it follows from \cite[Corollary~28.2.2]{rockafellar1970convex} and \cite[Corollary~28.3.1]{rockafellar1970convex} that $(\bar{y},\bar{z}, \bar{ w}) \in \cY\times \cZ\times \cW$ is an optimal solution to problem \eqref{eq-ADMM3-linear} if and only if there exists a Lagrange multiplier $\bar{x} \in \cX$ such that
\begin{equation}\label{genopt3B}
0\in \cF \bar{x} + \partial f(\bar{y}), \quad 0 \in \cG \bar{x} + \partial
g(\bar{z}),\quad \cH \bar{ x} -b = 0, \quad \cF^*\bar{y}+\cG^*\bar{z} +\cH^*\bar{w} -c=0.
\end{equation}
 Moreover, any $\bar{x}\in \cX$ satisfying \eqref{genopt3B} is an optimal solution to the dual problem \eqref{eq-ADMM3-linear-primal}.



We consider the following
sPADMM
for solving \eqref{eq-ADMM3-linear}.

\bigskip
\centerline{\fbox{\parbox{\textwidth}{
{\bf Algorithm sPADMM3c}: {\bf A convergent $3$-block semi-proximal ADMM for solving (\ref{eq-ADMM3-linear}).}
\\[5pt]
Let   $\sigma > 0$ and  $\tau \in (0,\infty)$ be given  parameters.
Choose $y^0 \in {\rm dom}(f)$, $z^0 \in {\rm dom}(g)$, and
$x^0 \in \cX$ such that $\cH  x^0= b$. Set $w^0 := (\cH \cH^* )^{-1} \cH (c-\cF^*  y^0-\cG^* z^0)$. Perform the $k$th iteration as follows:
\begin{description}
\item[Step 1.]  Compute
$\displaystyle y^{k+1} = \argmin  L_{\sigma}(y,z^k,w^k; x^k) +\frac{\sigma}{2}\norm{  y-y^k}^2 _{\cT_f}$.
\item[Step 2.]  Compute $\displaystyle w^{k+\frac{1}{2}}= \argmin  L_{\sigma}(y^{k+1},z^k,w; x^k)   = (\cH\cH^* )^{-1} \cH (c-\cF^*  y^{k+1}-\cG^* z^k)$ and
$\displaystyle z^{k+1} =\argmin  L_\sigma(y^{k+1},z,w^{k+\frac{1}{2}}; x^k) +\frac{\sigma}{2}\norm{ z-z^k}^2 _{\cT_g} $.
\item[Step 3.] Compute $w^{k+1}= \argmin  L_{\sigma}(y^{k+1},z^{k+1},w; x^k)   = (\cH \cH^* )^{-1} \cH (c-\cF^*  y^{k+1}-\cG^* z^{k+1})$.
     \item [Step 4.]  Compute  $x^{k+1} = x^k + \tau \sigma (\cF^* y^{k+1} + \cG^* z^{k+1} + \cH^* w^{k+1} - c)$.
\end{description}
}}}

\bigskip

Note that in Step 2 of Algorithm sPADMM3c, by direct calculations  we should have
 \[w^{k+\frac{1}{2}}= \argmin  L_{\sigma}(y^{k+1},z^k,w; x^k)  = (\cH \cH^* )^{-1} [\cH (c-\cF^*  y^{k+1}-\cG^* z^{k}) + \sigma^{-1}(b - \cH x^k )].
 \]
  However, by using Proposition \ref{prop:equivalence}, to be introduced later, we know that $ b - \cH x^k =0$ for all $k$. Thus
  $w^{k+\frac{1}{2}}= (\cH\cH^* )^{-1} \cH (c-\cF^*  y^{k+1}-\cG^* z^k)$.   In Step 3, $w^{k+1}$ is computed in a similar way.
When  $\cT_f=0$ and $\cT_g=0$, i.e., the proximal terms $\norm{y-y^k}_{\cT_f}^2$
and $\norm{z-z^k}_{\cT_g}^2$ are absent,
 Algorithm sPADMM3c  will become our convergent ADMM for solving (\ref{eq-ADMM3-linear}) (ADMM3c in short).  One reason for including $\cT_f$ and $\cT_g$ is to ensure that both $y^{k+1}$ and $z^{k+1}$ are well defined; see further discussions on this part in
   Section \ref{section:ADMM2B}. The difference between our ADMM3c and the directly extended 3-block ADMM (ADMM3d in short) is that
   we perform an extra intermediate step to compute $w^{k+\frac{1}{2}}$ before  computing $z^{k+1}$, i.e., at the $k$th iteration  we perform
   a particularly chosen   essentially BCD cycle in updating the variable $y, z, w$ in the terminology of Tseng \cite{tseng2001convergence}. Except for this extra step, ADMM3c is as simple as ADMM3d, which at each iteration
   performs a Gauss-Seidel BCD cycle in updating the variable $y, z, w$.
   Observe that in both ADMM3c and ADMM3d, we need to solve linear systems involving
the fixed operator $\cH\cH^*$.
For the case where the computation (which only needs to be done once) of $\cH\cH^*$ and
its (sparse) Cholesky factorization can be done at a moderate cost,
Step 3 of the above algorithm can be performed cheaply.
Now, under the condition that the Cholesky factorization of $\cH\cH^*$ is available,
the extra cost for computing $w^{k+\frac{1}{2}}$ is actually insignificant.
The reward for doing  the extra step in computing $w^{k+\frac{1}{2}}$
 is that  we are  able to  prove the convergence of our ADMM3c not only for $\tau =1$ but also allow $\tau$ to take a larger step-length, e.g,  $\tau =1.618$,
so as to achieve faster convergence than
 the directly extended ADMM3d.
Note that if the $z$-part disappears, then Step 2 of Algorithm sPADMM3c disappears and
our ADMM3c is identical  to  the classic 2-block ADMM but with $\tau \in (0,2)$   instead of  $\tau \in (0, ({1+\sqrt{5}}\,)/{2})$ due to   our requirement that $\cH x^0=b$.



Next we will prove the convergence of  Algorithm sPADMM3c for solving (\ref{eq-ADMM3-linear}) by relating it to the  generic 2-block
sPADMM
for solving a 2-block convex optimization problem discussed in Section \ref{section:ADMM2B}.

For problem (\ref{eq-ADMM3-linear}), one can obtain $w$ explicitly from the equality constraint $\cF^*   y + \cG^*  z + \cH^*  w = c$ as follows
\begin{eqnarray} \label{eq-y3}
w(y,z) = (\cH\cH^* )^{-1} \cH (c-\cF^*  y-\cG^* z), \quad (y,z) \in \cY \times \cZ.
\end{eqnarray}
Substituting (\ref{eq-y3})  into  (\ref{eq-ADMM3-linear}), we can recast  (\ref{eq-ADMM3-linear})  equivalently as
\begin{eqnarray}
   \min \left\{ f(y) + g(z) + \inprod{\bar{b}}{\cF^*  y+\cG^* z-c} \mid \cQ(\cF^*  y+\cG^* z-c) = 0 \right\},
  \label{eq-ADM3-no3}
\end{eqnarray}
where $\bar{b} := \cH^* (\cH \cH^* )^{-1}b$,
 \[
\cQ := \cI - \cP, \qquad \cP: = \cH^*  (\cH  \cH^* )^{-1} \cH
 \]
 and  $\cI:\cX \to \cX$ is the identity map.
It is easy to check that the two operators  $\cQ$ and $\cP$ satisfy the following
properties:
\begin{eqnarray} \label{eq:property-PQ}
\cP^* = \cP, \quad \cQ^* = \cQ, \quad \cP^*\cP = \cP, \quad \cQ^*\cQ = \cQ, \quad   \cP \cH^*= \cH^* ,\quad  \cH\cQ = 0,\quad  \cQ\cH^* = 0.
\end{eqnarray}
   The dual of (\ref{eq-ADM3-no3}) is given by
\begin{eqnarray} \label{eq-ADM3-primal1}
 \begin{array}{ll}
    \min & f^*(u)+g^*(v)+\inprod{c}{\bar{b}+\cQ\lambda}
       \\[5pt]
      {\rm s.t.} & \cF (\bar{b}+\cQ\lambda) +u =0, \quad
       \cG (\bar{b}+\cQ\lambda) + v =0.
 \end{array}
\end{eqnarray}
Note that (\ref{eq-ADM3-primal1}) is equivalent to (\ref{eq-ADMM3-linear-primal}) if we let $x=\bar{b}+\cQ\lambda$.

Let $\sigma>0$ be a positive constant. Define the augmented Lagrange function for (\ref{eq-ADM3-no3}) by
\begin{eqnarray}
 \hL_{\sigma}(y,z; \lambda) &=& f(y) + g(z) + \inprod{\bar{b}}{\cF^*  y+\cG^* z-c} \nonumber
 \\
 &&+\; \inprod{\lambda}{\cQ(\cF^*  y+\cG^* z-c)}
 + \frac{\sigma}{2} \norm{\cQ(\cF^*  y+\cG^* z-c)}^2.
\label{eq-ADM3-AL}
\end{eqnarray}
Now we can apply the generic
2-block
sPADMM
discussed in Section \ref{section:ADMM2B} to (\ref{eq-ADM3-no3}).

 \bigskip
\centerline{\fbox{\parbox{\textwidth}{
{\bf Algorithm sPADMM2s}: {\bf A specific 2-block semi-proximal ADMM for solving (\ref{eq-ADM3-no3}).}
\\[5pt]
Let   $\sigma > 0$ and  $\tau \in (0,\infty)$ be given parameters.
Choose $y^0 \in {\rm dom}(f)$, $z^0 \in {\rm dom}(g)$, and
$\lambda^0 \in {\rm Range}(\cQ)$.  Perform the $k$th iteration as follows:
\begin{description}
\item[Step 1.]  Compute
$\displaystyle y^{k+1} = \argmin \hL_{\sigma}(y,z^k; \lambda^k)+\frac{\sigma}{2}\norm{ \cF^*  (y-y^k)}^2_{\cP}+\frac{\sigma}{2}\norm{  y-y^k}^2 _{\cT_f}. $
\item[Step 2.]  Compute
$\displaystyle z^{k+1} = \argmin \hL_\sigma(y^{k+1},z; \lambda^k)+\frac{\sigma}{2}\norm{ \cG^* (z-z^k)}^2_{\cP} +\frac{\sigma}{2}\norm{  z-z^k}^2 _{\cT_g}.$
\item[Step 3.] Compute $\lambda^{k+1} = \lambda^k + \tau \sigma \cQ(\cF^* y^{k+1} + \cG^* z^{k+1}-c)$.
\end{description}
}}}
\bigskip

It is important to note that in Algorithm sPADMM2s  for solving (\ref{eq-ADM3-no3}),
we have two proximal terms in both Step 1 and Step 2 instead of one proximal term.  In particular, the first proximal term is necessary as neither $\cQ\cF^*$ nor $\cQ \cG^*$ is injective.  Next, we establish the equivalence of
 Algorithm sPADMM3c for solving problem (\ref{eq-ADMM3-linear}) and Algorithm sPADMM2s for solving problem (\ref{eq-ADM3-no3}).

 \begin{proposition}\label{prop:equivalence}
Let   $\sigma > 0$ and  $\tau \in (0,\infty)$ be given  parameters.
Choose $y^0 \in {\rm dom}(f)$, $z^0 \in {\rm dom}(g)$ and
$\lambda^0 \in {\rm Range}(\cQ)$.  Let $x^0=\cH^* (\cH \cH^* )^{-1}b+ \lambda^0$.  Then for any $k\ge 0$,  we have the following results
\begin{enumerate}
\item [{\rm (i)}]
the point $(y^{k}, z^{k})$ generated by Algorithm sPADMM2s for solving problem (\ref{eq-ADM3-no3}) is identical to the point
$(y^{k}, z^{k})$ generated by Algorithm  sPADMM3c for solving problem (\ref{eq-ADMM3-linear});

\item [{\rm (ii)}] $\lambda ^{k}$ and $x^{k}$ satisfy the relation
\[
x^{k}=\cH^* (\cH \cH^* )^{-1}b+ \lambda^{k}.
\]

\end{enumerate}
\end{proposition}

\begin{proof} We prove this proposition by induction.
First, note that since $\cH \cQ = 0$ and $\lambda^0 \in {\rm Range}(\cQ)$, we have $\cH x^0 = b$.

Recall that $w^k = (\cH\cH^*)^{-1}\cH(c-\cF^* y^k-\cG^* z^k)$.
Note that by direct computations we have $\cH^* w^k = \cP \cH^* w^k = \cP (c-\cF^* y^k-\cG^* z^k)$.
 Assume that the conclusions of this proposition hold up to $k\ge 0$.
Then, by using \eqref{eq:property-PQ} and the facts that $\lambda^k$ is in the range of $\cQ$ and $\cH x^k = b$, we can easily check  for every $k= 0,1, \ldots $ that
\begin{eqnarray*}
 && \hspace{-0.7cm}
y^{k+1} = \argmin \big\{ \hL_{\sigma}(y,z^k; \lambda^k)+\frac{\sigma}{2}\norm{ \cF^*  (y-y^k)}^2_{\cP}+\frac{\sigma}{2}\norm{  y-y^k}^2 _{\cT_f} \big\}
 \\[5pt]
 &=& \argmin \big\{ \hL_{\sigma}(y,z^k; \lambda^k)+\frac{\sigma}{2}\norm{ \cP(\cF^* y +\cG^*z^k +\cH^* w^k-c)}^2+\frac{\sigma}{2}\norm{  y-y^k}^2 _{\cT_f} \big\}
 \\[5pt]
  &=& \argmin \left\{ \begin{array}{l}
 f(y) + g(z^k) + \inprod{\bar{b}}{\cF^*  y+\cG^* z^k-c} +\; \inprod{\lambda^k}{\cQ(\cF^*  y+\cG^* z^k-c)}
 \\[3pt]
  + \frac{\sigma}{2} \norm{\cQ(\cF^*  y+\cG^* z^k-c)}^2 +\frac{\sigma}{2}\norm{ \cP(\cF^* y +\cG^*z^k +\cH^* w^k-c)}^2+\frac{\sigma}{2}\norm{  y-y^k}^2 _{\cT_f}
\end{array} \right\}
\\[5pt]
 &=& \argmin \left\{ \begin{array}{l}
  f(y) + g(z^k) + \inprod{x^k}{\cF^*  y+\cG^* z^k-c} + \frac{\sigma}{2} \norm{ \cQ(\cF^* y +\cG^*z^k +\cH^* w^k-c)}^2
 \\[3pt]
 +\frac{\sigma}{2}\norm{ \cP(\cF^* y +\cG^*z^k +\cH^* w^k-c)}^2+\frac{\sigma}{2}\norm{  y-y^k}^2 _{\cT_f}
\end{array}\right\}
\\[5pt]
 &=& \argmin \left\{ \begin{array}{l}
 f(y) + g(z^k) - \inprod{b}{w^k} + \inprod{x^k}{\cF^* y+\cG^* z^k+\cH^* w^k-c}
 \\[3pt]
+ \frac{\sigma}{2} \norm{ \cQ(\cF^* y +\cG^*z^k +\cH^* w^k-c)}^2
 +\frac{\sigma}{2}\norm{ \cP(\cF^* y +\cG^*z^k +\cH^* w^k-c)}^2+\frac{\sigma}{2}\norm{  y-y^k}^2 _{\cT_f}
\end{array}\right\}
\\[5pt]
 &=&  \argmin \left\{ \begin{array}{l}
f(y) + g(z^k) - \inprod{b}{w^k} + \inprod{x^k}{\cF^* y+\cG^* z^k+\cH^* w^k-c}
 \\[3pt]
+ \frac{\sigma}{2} \norm{ \cF^* y +\cG^*z^k +\cH^* w^k-c}^2+\frac{\sigma}{2}\norm{  y-y^k}^2 _{\cT_f}
\end{array}\right\}
\\[5pt]
&=& \argmin \big\{L_{\sigma}(y,z^k,w^k; x^k)  + \frac{\sigma}{2}\norm{  y-y^k}^2 _{\cT_f}\big\}.
\end{eqnarray*}
 Similarly, for every $k= 0,1, \ldots $, we have
\begin{eqnarray*}
 z^{k+1} &=&  \argmin \big\{\hL_\sigma(y^{k+1},z; \lambda^k)+\frac{\sigma}{2}\norm{ \cG^* (z-z^k)}^2_{\cP} +\frac{\sigma}{2}\norm{  z-z^k}^2 _{\cT_g}  \big\}
 \\[5pt]
 &=& \argmin  \big\{ \hL_{\sigma}(y^{k+1},z; \lambda^k)+\frac{\sigma}{2}\norm{ \cP(\cF^* y^{k+1} +\cG^*z +\cH^* w^{k+\frac{1}{2}}-c)}^2 +\frac{\sigma}{2}\norm{  z-z^k}^2 _{\cT_g}  \big\}
 \\[5pt]
&=& \argmin \big\{L_\sigma(y^{k+1},z,w^{k+\frac{1}{2}}; x^k)+\frac{\sigma}{2}\norm{  z-z^k}^2 _{\cT_g}\big\},
\end{eqnarray*}
where
$ w^{k+\frac{1}{2}}= \argmin  L_{\sigma}(y^{k+1},z^k,w; x^k)   = (\cH\cH^* )^{-1} \cH (c-\cF^*  y^{k+1}-\cG^* z^k)$
and
\begin{eqnarray*}
 \lambda^{k+1} &=& \lambda^k + \tau \sigma \cQ(\cF^*y^{k+1} + \cG^*z^{k+1}-c)
 \\[5pt]
 &=&x^k -\cH^* (\cH \cH^* )^{-1}b +  \tau \sigma (\cF^*y^{k+1} + \cG^*z^{k+1} +\cH^*w^{k+1}-c),
 \\[5pt]
  &=&x^{k+1} -\cH^* (\cH \cH^* )^{-1}b,
\end{eqnarray*}
 where
 $
 \displaystyle w^{k+1}= \argmin  L_{\sigma}(y^{k+1},z^{k+1},w; x^{k})   = (\cH\cH^* )^{-1} \cH (c-\cF^*  y^{k+1}-\cG^* z^{k+1}).
$
This completes our proof.
\end{proof}

Now we are ready to establish  the convergence results for Algorithm  sPADMM3c for solving (\ref{eq-ADMM3-linear}).

\begin{theorem}\label{thm:CPADMM3B}  Let
 $\Sigma_f$ and $\Sigma_g $    be   the  two self-adjoint and positive semidefinite operators defined by (\ref{monosub1}) and (\ref{monosub2}), respectively.
Suppose  that the solution set of problem (\ref{eq-ADMM3-linear}) is nonempty and that Assumption \ref{assumption:CQ1}  holds.   Assume that
$\cT_f$ and $\cT_g$ are chosen such
 that the  sequence $\{(y^k,z^k,w^k,x^k)\}$ generated by Algorithm sPADMM3c is well defined. Let $(\bar y, \bar z, \bar w)$ be any  optimal solution to (\ref{eq-ADMM3-linear}) and $\bar{x}$ be any optimal solution to \eqref{eq-ADMM3-linear-primal}, respectively.
For $k=0,1, 2, \ldots, $ denote
 \[
 y_e^k: =y^k -\bar y, \quad z_e^k: =z^k -\bar z, \quad w_e^k: =w^k -\bar w \quad {\rm and}\quad  x_e^k: = x^k -\bar {x}.
 \]
Then, under the condition either
(a)  $\tau\in (0, (1+\sqrt{5}\,)/2)$ or (b)  $\tau \ge (1+\sqrt{5}\,)/2$ but $\sum_{k=0}^\infty (\|\cG^*(z^{k+1}-z^k) +\cH^*(w^{k+1}-w^{k+\frac{1}{2}})\|^2 + \tau^{-1} \|\cF^* y^{k+1} +\cG^* z^{k+1}+\cH^* w^{k+1}-c\|^2) <\infty$,
the following results hold:
\begin{enumerate}
\item [{\rm (i)}] The sequence $\{\|x^{k+1}_e\|^2+\|z^{k+1}_e\|^2 _{(\sigma^{-1}\Sigma_g +\cT_g + \cG \cG^*)}+\|y^{k+1}_e\|^2 _{(\sigma^{-1}\Sigma_f +\cT_f +  \cF\cF^*)}\}$  is bounded.

 \item [{\rm (ii)}]   If $(y^\infty, z^\infty, w^\infty, x^\infty)$ is an accumulation point of $\{(y^k, z^k, w^k, x^k)\}$, then $(y^\infty, z^\infty, w^\infty)$ solves (\ref{eq-ADMM3-linear}) and $x^\infty$  solves \eqref{eq-ADMM3-linear-primal}, respectively and it holds that
     \[
     \lim_{k\rightarrow \infty}\left(  \|x^{k+1}_e\|^2+\|z^{k+1}_e\|^2 _{(\sigma^{-1}\Sigma_g +\cT_g + \cG \cG^*)}+\|y^{k+1}_e\|^2 _{(\sigma^{-1}\Sigma_f +\cT_f +  \cF \cF^*)}\right)=0,
     \]
     where in the definition of  $(y_e^k, z^k_e, w^k_e,x_e^k)$, the point $(\bar y, \bar z, \bar w, \bar x)$ is replaced by $(y^\infty, z^\infty,w^\infty, x^\infty)$.

   \item [{\rm (iii)}] If both $\sigma^{-1}\Sigma_f +\cT_f +  \cF\cF^*$ and $\sigma^{-1}\Sigma_g +\cT_g + \cG \cG^*$ are positive definite, then  the  sequence $ \{(y^k, z^k, w^k,x^k)\}$, which  is automatically well defined, converges to a unique limit, say, $(y^\infty, z^\infty, w^\infty, x^\infty)$  with $(y^\infty, z^\infty, w^\infty)$ solving (\ref{eq-ADMM3-linear}) and $x^\infty$ solving \eqref{eq-ADMM3-linear-primal}, respectively.
 \item [{\rm (iv)}]  When the the $z$-part disappears, the corresponding results in parts {\rm (i)--(iii)} hold for
 any $\tau \in (0,2)$ or $\tau \ge 2$ but $\sum_{k=0}^\infty \|\cF^* y^{k+1} + \cH^* w^{k+1} -c\|^2<\infty$.
\end{enumerate}
\end{theorem}
\begin{proof} By combing Theorem  \ref{thempadm} with Proposition \ref{prop:equivalence} and using  the relation (\ref{eq-y3}), we can readily obtain the conclusions of this theorem.
\end{proof}

\begin{remark}
     The main idea for proving the convergence of Algorithm sPADMM3c is via showing  that Algorithm sPADMM3c
      is equivalent to  Algorithm sPADMM2s, which is obtained by applying  Algorithm sPADMM2
      to the   problem \eqref{eq-ADM3-no3} using two special  semi-proximal terms $\cS = \cF\cP\cF^*+\cT_f$ and $\cT = \cG\cP\cG^*+\cT_g$ in Step 1 and Step 2, respectively. This simple discovery of the equivalence of Algorithm sPADMM3c for solving   (\ref{eq-ADMM3-linear}) and Algorithm sPADMM2s for solving its equivalent problem (\ref{eq-ADM3-no3}) is significant since  on the one hand it settles the convergence  of Algorithm sPADMM3c    by using known convergence results  for Algorithm sPADMM2s and on the other hand, it allows one to take   advantage of  the extremely  simple structure of   Algorithm sPADMM3c   in searching for an efficient implementation for solving convex conic programming.
      Note that one cannot even prove the convergence of Algorithm ADMM3c (without the two semi-proximal terms $\cT_f$ and $\cT_g$) by directly applying the classic $2$-block ADMM to \eqref{eq-ADM3-no3}. Actually, since neither $\cQ\cF^*$
      nor  $\cQ\cG^*$ is injective, one cannot use  the classic $2$-block ADMM to solve \eqref{eq-ADM3-no3} at all unless additional conditions on $f$ and $g$ are imposed.
\end{remark}

\subsection{Extensions}\label{section:extension}

In this subsection, we first consider the  $3$-block convex optimization problem (\ref{eq-ADMM3-linear}), i.e.,
  \begin{eqnarray}
 \min \left\{ f(y)+g(z)+h({w}) \mid \cF^* y + \cG^* z + \cH^* w = c \right\}
  \nonumber
\end{eqnarray}
for the case where at least  one of the three functions $f$, $g$ and $h$  is partially linear.
Without loss of generality, we assume that $h:\cW \equiv \cW_I \times \cW_E \to (-\infty, \infty]$ is of  the following partial linear structure
\[
{h(w)} =\theta(w_I)-\inprod{b}{w} =\theta(w_I)-\inprod{b_I}{w_I}-\inprod{b_E}{w_E}  \quad \forall \, w \equiv (w_I, w_E) \in \cW_I \times \cW_E,
\]
where $b\equiv (b_I,b_E)\in \cW_I\times \cW_E$ is a given vector and $\theta: \cW_I\to (-\infty, \infty]$ is a closed proper nonlinear convex function. Decompose
$\cH \equiv \cH_I\times \cH_E$ such that  for any $x\in \cX$, \[
 \left(\begin{array}{c}\cH_I x \\ \cH_E x \end{array}\right)
\equiv  \cH x
= \cH_I x \times \cH_E x \in   \cW_I \times  \cW_E.
\]
Again, by removing redundancy in $\cH$ if necessary, we   assume that $\cH_E\cH_E^*$ is invertible, i.e., $\cH_E$ is surjective.
We also assume that the
Cholesky factorization of $\cH_E\cH_E^*$ can be computed at a moderate cost. In order to apply our proposed convergent  3-block semi-proximal ADMM to  solve
  \begin{eqnarray}
 \min \left\{ f(y)+g(z)+h({w}) \mid \cF^* y + \cG^* z + \cH_I^* w_I +\cH_E^*w_E = c \right\},
  \label{eq:partial-linear-h}
\end{eqnarray}
we need to convert it into  the form of problem (\ref{eq-ADMM3-linear}) with $h$ being linear.  For this purpose,  we define $\cH_3: \cW\to \cX \times \cW_I $  to be the following linear map
\[
\cH_3^* w: = \left(\begin{array}{l}  \cH_I^*w_I +\cH_E^* w_E\\ -\cD_I^* w_I\end{array} \right) \quad \forall\, w\in \cW,
\]
where   $\cD_I: \cW_I \to \cW_I$ is a particularly  chosen nonsingular linear operator, e.g., the identity operator $\cI:\cW_I \to \cW_I$,  and $\cD_I^*$ is the adjoint of $\cD_I$.

 We consider the following two cases.

\medskip
{\bf Case 1)}.  The inverse of $\cH_3\cH_3^*$  can be computed at a moderate cost.

 By introducing a slack variable $\eta\in \cW_I$, we  can rewrite
(\ref{eq:partial-linear-h}) as
  \begin{eqnarray}
 \min \left\{f(y)+g(z) + \theta(w_I)- \inprod{b}{w} \mid \cF^* y + \cG^* z + \cH_I^*w_I +\cH_E^* w_E = c, \;  \cD_I^*( \eta-w_I) =0\right\},
\label{eq-ADMM3-A0}
\end{eqnarray}
which can then be cast into the form of problem (\ref{eq-ADMM3-linear}) as follows:
\begin{eqnarray}
 \min \left\{\Big(f(y)+\theta(\eta)\Big) +g(z) - \inprod{b}{w}  \mid \left (\begin{array}{c} \cF^* y \\ \cD_I^* \eta\end{array}\right)  + \left (\begin{array}{c} \cG^* z\\ 0\end{array}\right) + \cH_3^*w= \left (\begin{array}{c} c\\ 0\end{array}\right) \; \right\}.
 \label{eq-ADMM3-A}
\end{eqnarray}
The convergent 3-block
sPADMM discussed in Section \ref{subsection:h-linear} for solving  problem (\ref{eq-ADMM3-linear}) can then be applied to problem  (\ref{eq-ADMM3-A}) in a straightforward way.

\medskip
{\bf Case 2)}. The inverse of $\cH_3\cH_3^*$  cannot be computed at a moderate cost.

Let $\cD: \cX \to \cX$ be a given nonsingular linear  operator and $\cD^*$ be its adjoint. We assume that $\cD$ is chosen in such a way that the inverse of $\cI + \cD \cD^*$ can be computed explicitly with low costs (e.g., $\cD= 5\cI$).   By introducing a slack variable $s\in \cX$, we  can rewrite  \eqref{eq:partial-linear-h}
 as
  \begin{eqnarray}
 \min \left\{f(y)+g(z) + \theta(w_I)- \inprod{b}{w} \mid \cF^* y + \cG^* z + s +\cH_E^* w_E = c, \;  \cD^*( \cH_I^*w_I-s) =0\right\},
\label{eq-ADMM3-A0}
\end{eqnarray}
which can then  be recast in the form of (\ref{eq-ADMM3-linear}) as follows:
\begin{eqnarray}
 \min \left\{\Big(f(y) + \tilde \theta(w_I) \Big)+g(z)- \inprod{b_E}{w_E}  \mid \left (\begin{array}{c} \cF^* y \\ \cD^* \cH_I^*w_I\end{array}\right)  + \left (\begin{array}{c} \cG^* z\\ 0\end{array}\right) + \cB^* (s,w_E)= \left (\begin{array}{c} c\\ 0\end{array}\right) \; \right\},
 \label{eq-ADMM3-B}
\end{eqnarray}
where the convex function ${\tilde \theta}(\cdot)\equiv \theta(\cdot) -\inprod{b_I}{\cdot}$ and the linear map $\cB: \cX \times \cW_E\to \cX \times \cX$  is defined by
\[
\cB^* (s,w_E) := \left(\begin{array}{l}  s +\cH_E^* w_E\\ -\cD^*s\end{array} \right) \quad \forall\, (s, w_E)\in \cX\times \cW_E.
\]
As in Case 1), we can apply the convergent 3-block
sPADMM discussed in Section \ref{subsection:h-linear}   to problem  (\ref{eq-ADMM3-B}) in a straightforward way as now the inverse of $\cB\cB^*$ can be computed based on $(\cI + \cD\cD^*)^{-1}$
and the inverse of  $\cH_E(\cI- (\cI + \cD\cD^*)^{-1}) \cH_E^*$.
In our numerical experiments in Section \ref{section:numericalresults},
we choose $\cD =\alpha \cI$ for some $\alpha \in [3,6]$, and $\alpha$
is dynamically adjusted according to the progress of the
algorithm.

Though not the focus of this paper,
here we will also briefly explain   how to
extend our  convergent semi-proximal ADMM to deal with the general $q$-block convex optimization
 problem (\ref{eq-ADMMq}):
  \begin{eqnarray}
 \min \left\{ {\sum}_{i=1}^q\theta_i(w_i)  \mid {\sum}_{i=1}^q\cB^*_i w_i = c \right\}.
  \nonumber
\end{eqnarray}

For any $i\ge 3$, let $\cD_i: \cX \to \cX$  be a given nonsingular linear  operator and $\cD_i^*$ be its adjoint.    By introducing slack variables $s_i\in \cX$, $i=3, \ldots, q$, we  can then rewrite
(\ref{eq-ADMMq}) equivalently  as
 \begin{eqnarray}
 \min \left\{ {\sum}_{i=1}^q\theta_i(w_i)  \mid \cB_1^* w_1 +\cB_2^* w_2+  s_3+\cdots +s_q = c, \; \cD^*_i( \cB_i^*w_i - s_i) = 0, \ i=3, \ldots, q \right\}.
 \label{eq-ADMMq-A}
\end{eqnarray}
Let $l: =\lfloor q/2\rfloor$ be the largest integer that is smaller or equal to $q/2$. Define
$h(s_3, \ldots, s_q)\equiv 0$,
\begin{eqnarray*}
&f(w_1,w_3 \ldots, w_{l+1}) = \theta_1(w_1)+\sum_{i=3}^{l+1}\theta_i(w_i), \quad
g(w_2,w_{l+2}, \ldots, w_q)=\theta(w_2)+\sum_{i=l+2}^q\theta_i(w_i), &
\\[5pt]
&	\cF^*(w_1,w_3,\dots,w_{l+1}) =\left(\begin{array}{c}
		\cB_1^*w_1 \\
		\cD_3^*\cB_3^* w_3
		\\
		\vdots
		\\
		\cD_{l+1}^*\cB_{l+1}^* w_{l+1}
		\\
		0
		\\
		\vdots
		\\
		0
	\end{array} \right), \;
	\cG^*(w_2,w_{l+2},\dots,w_{q}) =\left(\begin{array}{c}
		\cB_2^*w_2 \\
		0
		\\
		\vdots
		\\
		0
		\\
		\cD_{l+2}^*\cB_{l+2}^*w_{l+2}
		\\
		\vdots
		\\
		\cD_{q}^*\cB_{q}^* w_{q}
	\end{array} \right) &
\end{eqnarray*}
 for any $w_i\in \cW_i$, $i=1, \ldots, q$ and $s_j\in \cX$, $j=3, \ldots, q$. Then it is easy to see that (\ref{eq-ADMMq-A}) is in the form of (\ref{eq-ADMM3-linear}) with
 \begin{eqnarray*}
 \cH^*(s_3, \ldots, s_q) =\left(\begin{array}{c}s_3 +\cdots +s_q \\
 -\cD_3^*  s_3
 \\
 \vdots
 \\
 -\cD_q^*s_q
 \end{array} \right).
 \end{eqnarray*}
Note that we have
\begin{eqnarray*}
 \cH\cH^*  = \left( \begin{array}{ccc}
\cD_3\cD_3^* & & \\
 & \ddots & \\
 & & \cD_q\cD_q^*
\end{array}\right) + \left( \begin{array}{c} \cI \\ \vdots \\ \cI \end{array}\right)
 \left( \begin{array}{c} \cI \\ \vdots \\ \cI \end{array}\right)^*.
\end{eqnarray*}
Thus, if for each $i\in \{3,\dots,q\}$, $\cD_i$ is chosen such that $\cD_i\cD_i^*$ can be inverted easily, e.g.,
 $\D_i =\alpha \cI$ for some $\alpha>0$, then we can compute the inverse of $\cH\cH^*$ analytically via the Sherman-Morrison-Woodbury formula if
$\cI+\sum_{i=3}^q (\cD_i\cD_i^*)^{-1}$ is also easy to invert. In this way,
 our convergent 3-block  semi-proximal ADMM discussed in Section \ref{subsection:h-linear}  can then be applied to problem \eqref{eq-ADMMq-A}
 directly.

\section{Applications to   conic programming}\label{section:ConicProgramming}

In this section, we show how to apply  our convergent $3$-block
sPADMM to solve  conic programming ({\bf P}) and
its dual   ({\bf D}).
Here we always assume that $\cA_E\cA_E^*$ is invertible and
its Cholesky factorization can be computed at a moderate cost.

\subsection{The case where $m_I=0$ }\label{subsection:ConicD-3block}

In this subsection, we show how our convergent ADMM3c can be used to solve   conic programming ({\bf P}) without the inequality constraints $\cA_I x_I \geq b_I$. In this case,  the conic programming ({\bf P}) reduces to
\begin{eqnarray}
  \max \left\{ \inprod{-c}{x} \mid \cA_E x=  b_{ E}, \;   x \in \cK, \;   x\in \cK_p \right\}
  \label{eq-primal3}
\end{eqnarray}
with its dual taking the form of
\begin{eqnarray}
 \quad  \min \left\{\delta_{\cK^*}(s) + \delta_{\cK_p^*}(z)  - \inprod{b_E}{y_E}\mid s +z+ \cA^*_E y_E = c\right\}.
  \label{eq-dual3}
\end{eqnarray}
For   problem \eqref{eq-dual3}, instead of using  the constraint qualification  imposed in Assumption \ref{assumption:CQ1}, we will use the following more familiar one in the conic programming field.

\begin{assumption} \label{assumption-feasible3block}
{\rm (a)} For problem (\ref{eq-primal3}), there exists a feasible solution $\hat{x} \in \cK$  such that
\begin{eqnarray*}
\mathcal{A}_E\hat{x} = b_E, \; \hat{x}\in {\rm int}(\cK), \; \hat{x} \in \cK_p \,.
\end{eqnarray*}
{\rm (b)} For problem (\ref{eq-dual3}), there exists a feasible solution $(\hat{s}, \hat{z}, \hat{y}_E) \in  \cK^* \times \cX\times \Re^{m_E}$ such that
\begin{eqnarray*}
\hat{s} + \hat{z}+ \cA^*_E{\hat{y}_E}  = c, \; \hat{s}\in {\rm int}(\cK^*), \; \hat{z} \in \cK_p^* \,.
\end{eqnarray*}
\end{assumption}

It is known from convex analysis (e.g, \cite[Corollary 5.3.6]{borwein2006convex}) that under Assumption \ref{assumption-feasible3block}, the strong duality for  \eqref{eq-primal3} and  \eqref{eq-dual3}
holds and the following Karush-Kuhn-Tuck  (KKT) condition has nonempty solutions:
\begin{eqnarray} \label{eq-KKT}
\left\{ \begin{array}{l l}
  \cA_E x -b_E =  0, \quad   s + z +\cA_E^*y_E  - c = 0, &    \\[5pt]
   \inprod{x}{s}  =  0, \; x\in \cK,\; s\in \cK^*, \;
   \inprod{x}{z}  =  0, \; x\in \cK_p,\; z\in \cK_p^* .  &
 \end{array}
 \right.
\end{eqnarray}

Let $\sigma > 0$ be given. The augmented Lagrange
 function for (\ref{eq-dual3}) is defined by
\begin{eqnarray}
L_\sigma(s,z,y_E;x): &=& \delta_{\cK^*} (s) +  \delta_{\cK^*_p} (z)+ \inprod{-b_E}{y_E} + \inprod{x}{s + z + \cA_E^*y_E - c} \nonumber
\\
&&+ \frac{\sigma}{2} \norm{s + z + \cA_E^*y_E -c}^2, \quad (s,z,y_E,x) \in \cX \times \cX \times \Re^{m_E}\times \cX.
\label{eq-Conic-ADM-3-AL}
\end{eqnarray}

We can apply our convergent 3-block ADMM (without the proximal terms) to  problem (\ref{eq-dual3}) to obtain the following algorithm.

\bigskip
\centerline{\fbox{\parbox{\textwidth}{
{\bf Algorithm Conic-ADMM3c}: {\bf A convergent  3-block ADMM for solving  (\ref{eq-dual3}).}
\\[5pt]
Given   parameters $\sigma > 0$ and $\tau \in (0,\infty)$.
Choose
$s^0\in \cK^*$, $z^0\in \cK_p^*$, $x^0 \in \cX$ such that $\cA_E x^0 = b_E$. Set $y_E^{0} = (\cA_E\cA_E^*)^{-1} \cA_E(c-s^0-z^0)$.  Perform the $k$th iteration as follows:
\begin{description}
\item[Step 1.]  Compute $\displaystyle s^{k+1}= \argmin  L_{\sigma}(s,z^k, y_E^k; x^k)
  = \Pi_{\cK^*}\left(c-z^k-\cA_E^* y_E^{k} - \sigma^{-1}x^k \right)$.
\item[Step 2.]  Compute $y_E^{k+\frac{1}{2}} = \argmin  L_{\sigma}(s^{k+1},z^k, y_E; x^k)   =(\cA_E\cA_E^*)^{-1} \cA_E(c-s^{k+1}-z^k)$ and
\[ z^{k+1}= \argmin  L_{\sigma}(s^{k+1},z, y_E^{k+\frac{1}{2}}; x^k)    = \Pi_{\cK^*_p}\left(c-s^{k+1}-\cA_E^* y_E^{k+\frac{1}{2}} - \sigma^{-1}x^k \right).\]
\item[Step 3.]  Compute $y_E^{k+1} = \argmin  L_{\sigma}(s^{k+1},z^{k+1}, y_E; x^k)   =(\cA_E\cA_E^*)^{-1} \cA_E (c-s^{k+1}-z^{k+1})$.
 \item[Step 4.] Compute  $x^{k+1} = x^k+\tau\sigma(s^{k+1}  + z^{k+1}+\cA_E^*y_E^{k+1} -c)$.
\end{description}
}}}

\medskip
The following convergence results for  Algorithm Conic-ADMM3c  for solving problem  (\ref{eq-dual3})  can be derived directly from Theorem \ref{thm:CPADMM3B}.

\begin{theorem}\label{thm:ADMM3c}
Suppose  that   Assumption \ref{assumption-feasible3block} holds and that $\cA_E$ is surjective.  Then the sequence $ \{(s^k, z^k, y_E^k,x^k)\}$ generated by Algorithm Conic-ADMM3c is well defined.  Furthermore, under the condition that either
(a)  $\tau\in (0, (1+\sqrt{5}\,)/2)$ or (b)  $\tau \ge (1+\sqrt{5}\,)/2$ but $\sum_{k=0}^\infty (\|(z^{k+1}-z^k) +\cA_E^*(y_E^{k+1}-y_E^{k+\frac{1}{2}})\|^2 + \tau^{-1} \|s^{k+1} +z^{k+1}+\cA_E^* y_E^{k+1}-c\|^2) <\infty$,
  it  holds that
\begin{enumerate}
\item[{\rm (i)}] The sequence $ \{(s^k, z^k, y_E^k,x^k)\}$  converges to a unique limit, say, $(s^\infty, z^\infty, y_E^\infty, x^\infty)$ satisfying the KKT condition (\ref{eq-KKT}).
 \item [{\rm (ii)}]  When   $\cK_p =\cX$, i.e., the $z$-part disappears, the corresponding result in part {\rm (i)} holds for any $\tau \in (0,2)$ or $\tau \ge 2$ but $\sum_{k=0}^\infty \|s^{k+1} + \cA_E^* y_E^{k+1} -c\|^2<\infty$ .
\end{enumerate}
\end{theorem}
\subsection{The case where $m_I>0$ }\label{subsection:ConicD-4block}

Here we consider conic programming problem ({\bf P}) with  inequality constraints, i.e.,
\begin{eqnarray}
  \max \left\{ \inprod{-c}{x} \mid \cA_E x=  b_{ E}, \; \cA_I x \ge b_I ,\; x \in \cK, \;   x\in \cK_p \right\}.
  \label{P:primal}
\end{eqnarray}
If $m_I$, the number of inequality constraints  $\cA_I x \ge b_I$, is relatively small,  we can introduce a slack variable  to  convert (\ref{P:primal}) into the form of problem (\ref{eq-primal3}) with three blocks of constraints and then apply Algorithm Conic-ADMM3c introduced in Section \ref{subsection:ConicD-3block} to solve it. We omit the details here.

Next,  we consider the case where $m_I$ is large. The dual of (\ref{P:primal}) is given by
\begin{eqnarray}
  \min \left\{ \delta_{\cK^*} (s) +  \delta_{\Re^{m_I}_+} (y_I)+ \delta_{\cK^*_p} (z)-\inprod{b}{y} \mid s +  \cA_I^*y_I + z +\cA_E^*y_E = c\;
    \right\}.
  \label{eq-dual4}
\end{eqnarray}
Let $\cD:\cX \to \cX$ be a   given nonsingular linear operator    and $\cD^*$ be its adjoint.
 In this case, we can rewrite (\ref{eq-dual4}) equivalently as
 \begin{eqnarray}
  \min \left\{ \delta_{\cK^*} (s) +  \delta_{\Re^{m_I}_+} (y_I)+ \delta_{\cK^*_p} (u) -\inprod{b}{y} \mid s +  \cA_I^*y_I + z +\cA_E^*y_E = c, \; \cD^*(u-z)=0\;
    \right\}.
  \label{eq-dual4-X}
\end{eqnarray}
Define ${ \cB}:\cX\times \Re^{m_E} \rightarrow \cX\times \cX$ to be  the linear
map whose adjoint ${ \cB}^*$ satisfies
\[
{ \cB}^*(z,y_E)=\left (\begin{array}{l}  z+ { \cA}_E^*y_E \\
  - \cD^* z \end{array}\right).
\]
Then problem (\ref{eq-dual4-X}) can be reformulated as
\begin{eqnarray}
 \min \left\{ f(s,u) + g(y_I) + h(z, y_E)
 \mid
 \left (\begin{array}{c} s \\ \cD^* u\end{array}\right)+\left (\begin{array}{c} \cA_I^*y_I \\ 0 \end{array}\right) +\cB^*(z,y_E) =\left (\begin{array}{c} c \\ 0\end{array}\right)
     \right\},
  \label{eq-dual01-AoA}
\end{eqnarray}
 where
 \[
 f(s,u):=\delta _{{\cK}^*} (s)   + \delta _{\cK_p^*} (u) , \quad g(y_I): = \delta _{\Re_+^{m_I}} (y_I) -  \inprod{b_I}{y_I}, \quad
 h(z, y_E) := -  \inprod{b_E}{y_E}
 \]
for any $(s, u) \in \cX \times \cX$, $y_I \in \Re^{m_I}$ and $(z, y_E)\in \cX\times \Re^{m_E}$.
Note that since
\[
\cB\cB^* = \left(\begin{array}{cc} \cI +\cD\cD^* & \cA_E^* \\ \cA_E & \cA_E\cA^*_E
\end{array}\right)
\]
and the inverse of $\cA_E\cA_E^*$ is assumed to be computable at a moderate cost,
the inverse of  ${\cB}{\cB}^*$ can also be computed based on
$(\cI+\cD\cD^*)^{-1}$ and the inverse of
$\cA_E(\cI -  (\cI +\cD\cD^*)^{-1})\cA_E^*$.
For example, if $\cD$ is a simple nonsingular matrix (e.g., $\cD = \alpha \cI$ for some $\alpha>0$), then
the inverse of  ${\cB}{\cB}^*$ can be computed at a low cost once the inverse of $\cA_E\cA_E^*$ is available. Let $\rho_{\max}$ be the largest eigenvalue of the self-adjoint positive semidefinite operator $A_I\cA_I^*$. Then  we can apply our convergent Algorithm sPADMM3c given in Section \ref{subsection:h-linear}  directly  to problem  (\ref{eq-dual01-AoA})   by
defining
\[
\cT_f \equiv 0 \quad {\rm and} \quad \cT_g \equiv \rho_{\max} \cI -\A_I\cA_I^*
\]
to obtain a convergent $3$-block
sPADMM, denoted by Conic-sPADMM3c,  for solving   conic programming ({\bf P}).
The motivation for choosing the specific positive semidefinite linear operator $\cT_g$ above
is to make the computation of $y_I^{k+1}$ simple. As before,
the convergence analysis for
Algorithm Conic-sPADMM3c can be analyzed similarly as in Theorem \ref{thm:ADMM3c}. For  simplicity,  we omit the details here.

\section{Numerical experiments for SDP}\label{section:numericalresults}

In the last section, we have shown how our proposed   convergent $3$-block
sPADMM for solving the convex optimization problem \eqref{eq-ADMM3-linear} can be used to solve ({\bf P}) with/without the inequality constraints. In this section, we use ({\bf SDP}) problems as our test examples. We separate our test examples  into two groups. The  first group is for DNN-SDP without the inequality constraints $\cA_I x_I \geq b_I$. In the second group, we consider SDP problems arising  from 
relaxation of
binary integer quadratic (BIQ) programming  problems with a large number of inequality constraints $\cA_I x_I \geq b_I$.

\subsection{Numerical results for doubly non-negative SDP}

The doubly non-negative SDP takes the form of
\begin{eqnarray}
{\bf (DNN{\texttt -}SDP)} \quad  \max \left\{ \inprod{-C}{X} \mid\cA_E X=  b_{ E}, \; X \in \cS_+^n, \; X\in \cK_p \right\},
  \label{eq-primal-DNNSDP}
\end{eqnarray}
whose dual can be written as
\begin{eqnarray}
 \quad  \min \left\{\delta_{\cS^n_+}(S) + \delta_{\cK_p^*}(Z)  - \inprod{b_E}{y_E}\mid S +Z+ \cA^*_E y_E = C\right\}.
  \label{eq-dual-DNNSDP}
\end{eqnarray}
Obviously,  our proposed Algorithm Conic-ADMM3c given in Section \ref{subsection:ConicD-3block} for     solving  problem (\ref{eq-dual3}) can be applied to problem (\ref{eq-dual-DNNSDP}) directly.

\subsubsection{Doubly non-negative SDP problem sets}

In our numerical experiments, we test the following classes of doubly non-negative SDP problems.

\bigskip
(i) DNN-SDP problems arising from
the relaxation of a binary integer nonconvex quadratic (BIQ) programming:
\begin{eqnarray}
\min\Big\{ \frac{1}{2}x^TQx + \inprod{c}{x} \mid  x\in\{0,1\}^{n-1}\Big\}.
\end{eqnarray}
This problem has  been shown in \cite{burer2009copositive}  that under some mild
assumptions, it can equivalently be reformulated as the following
completely positive programming (CPP) problem:
\begin{eqnarray}
\min\Big\{ \frac{1}{2} \inprod{Q}{Y} + \inprod{c}{x} \mid
 \diag(Y) = x, X = [Y, x; x^T, 1]\in
{\cal C}_{pp}^n\Big\},
\end{eqnarray}where ${\cal C}_{pp}^n$ denotes the $n$-dimensional completely positive cone.
It is well known that even though ${\cal C}_{pp}^n$ is convex, it is
computationally intractable. To solve the CPP problem,  one would typically
relax ${\cal C}_{pp}^n$ to $\S^n_+\cap \cK_{\cP}$, and the relaxed problem
has the form of ({\bf SDP}):
\begin{eqnarray}
\begin{array}{ll}
   \min & \frac{1}{2} \inprod{Q}{Y} + \inprod{c}{x} \\[5pt]
  {\rm s.t.} &  {\rm diag}(Y) - x = 0,\quad \alpha = 1,\quad
   X = \left[ \begin{array}{cc} Y & x \\ x^T & \alpha  \end{array}\right]\in \S^n_+,\quad
  X \in \cK_{\cP},
  \end{array}
  \label{eq-BIQ-SDP}
\end{eqnarray}
where the polyhedral cone $\cK_{\cP} = \{X \in \cS^n \mid X \geq 0 \}$.
In our numerical experiments, the test data for $Q$ and $c$ are taken from Biq Mac Library  maintained by Wiegele, which is available at {\tt http://biqmac.uni-klu.ac.at/biqmaclib.html}

\bigskip
(ii) DNN-SDP problems arising  from the relaxation of maximum stable set problems. Given a graph $G$ with edge set $\mathcal{E}$, the SDP relaxation $\theta_+(G)$ of the maximum stable set problem is given by
\begin{eqnarray}
 \theta_+(G)  = \max \{ \inprod{ee^T}{X} \,\mid \ \inprod{E_{ij}}{X} =0,
  (i,j) \in \mathcal{E}, \,\inprod{I}{X} =1, \, X \in \S^n_+,\, X \in \cK_{\cP} \},
  \label{eq-theta-plus}
\end{eqnarray}
where $E_{ij}=e_ie_j^{\T} + e_je_i^{\T}$ and $e_i$ denotes the $i$th column
 of the $n\times n$ identity matrix, and $\cK_{\cP} = \{X \in \cS^n \mid X \geq 0 \}$.
In our numerical experiments, we test the graph instances $G$ considered in
\cite{Sloane}, \cite{toh2004}, and \cite{DIMACS}.

\bigskip
(iii) DNN-SDP problems arising from computing lower
bounds for quadratic assignment problems (QAPs). Let $\Pi$ be the set of $n\times n$ permutation matrices. Given
matrices $A, B \in \S^n$, the quadratic assignment
problem is given by
\begin{eqnarray}
 \bar v_{\rm QAP} := \min \{ \langle X, AXB \rangle \, :\, X \in \Pi\}.
 \label{eq-QAP}
\end{eqnarray}
For a matrix $X = [x_1,\dots,x_n]\in \Re^{n\times n}$, we will
identify it with the $n^2$-vector $x = [x_1; \dots; x_n]$. For a
matrix $Y\in R^{n^2 \times n^2}$, we let $Y^{ij}$ be the $n\times n$
block corresponding to $x_i x_j^T$ in the $n^2\times n^2$ matrix $xx^T$. It is shown
in \cite{povh2009copositive} that $\bar v_{\rm QAP}$ is bounded below by the
following number generated from the SDP relaxation of (\ref{eq-QAP}):
\begin{eqnarray}
  \begin{array}{rl}
  v \,:=\, \min &\displaystyle \inprod{B \otimes A}{Y} \\[5pt]
  {\rm s.t.} &  \sum_{i=1}^n Y^{ii} = I, \quad
  \inprod{I}{Y^{ij}} = \delta_{ij} \quad \forall \, 1\leq i \leq j\leq n, \\[5pt]
  &
  \inprod{E}{Y^{ij}} = 1\quad \forall\, 1\leq i \leq j\leq n, \\[5pt]
  & Y\in \S^{n^2}_+,\; Y \in \cK_{\cP},
  \end{array}
  \label{eq-PR}
\end{eqnarray}
where $E$ is the matrix of ones, and $\delta_{ij} = 1$ if $i=j$,
and $0$ otherwise, $\cK_{\cP} = \{X \in \cS^{n^2} \mid X \geq 0 \}$.
In our numerical experiments, the test instances $(A,B)$ are taken
from the QAP Library \cite{QAPLIB}.

\bigskip
(iv) DNN-SDP relaxation of clustering problems (RCPs)  described in
\cite[eq. (13)]{peng2007approximating}:
\begin{eqnarray}
\min \Big\{ \inprod{W}{I-X} \mid Xe = e, \inprod{I}{X} = K, X \in \S^n_+,
X \in \cK_{\cP} \Big\},
\end{eqnarray}
where $W$ is the so-called affinity matrix whose entries represent the
pairwise  similarities of the objects in the dataset, $e$ is the vector of ones, and $K$ is the number of clusters, $\cK_{\cP} = \{X \in \cS^n \mid X \geq 0 \}$.  All the data sets we test are from the UCI Machine Learning Repository (available at {\tt http://archive.ics.uci.edu/ml/datasets.html}). For some large size data sets, we only select the first $n$ rows. For example, the original data set ``{\tt spambase}" has 4601 rows and we select the first 1500 rows to obtain the test problem ``{\tt spambase-large.2}"
for which the number ``{\tt 2}" means that there are $K = 2$ clusters.

\bigskip
(v) DNN-SDP problems arising from semidefinite relaxation of frequency assignment problems (FAPs) \cite{eisenblatter2002frequency}.  Given a network represented by a graph $G$ and an edge-weight matrix $W$, a certain type of frequency assignment problem on $G$  can be relaxed into the following SDP (see \cite[eq. (5)]{burer2003computational}):
\begin{eqnarray} \label{eq-fap}
\begin{array}{ll}
   \max & \inprod{(\frac{k-1}{2k})L(G,W)-\frac{1}{2}{\rm Diag}(We)}{ X} \\[5pt]
  {\rm s.t.} & {\rm diag}(X) = e, \quad X \in \S^n_+,
  \\[3pt]
  & \inprod{-E^{ij}}{ X} = 2/(k-1) \quad \forall (i,j) \in U \subseteq E,
  \\[3pt]
  & \inprod{-E^{ij}}{ X} \leq 2/(k-1) \quad \forall (i,j) \in E \setminus U,
  \end{array}
\end{eqnarray}
where $k>1$ is an integer, $L(G,W) ;= {\rm Diag}(We)-W$ is the Laplacian matrix, $E^{ij}=e_ie_j^T + e_je_i^T$ with $e_i \in \Re^n$ the vector with all zeros except in the $i$th position, and $e \in \Re^n$ is the vector of ones.
Denote
\begin{eqnarray*}
M_{ij} =
\begin{cases}
-\frac{1}{k-1} &  \forall (i,j) \in E, \\
0 & {\rm otherwise}.
\end{cases}
\end{eqnarray*}
Then (\ref{eq-fap}) is equivalent to
\begin{eqnarray} \label{eq-fap-sdpp}
 \begin{array}{ll}
    \max & \inprod{(\frac{k-1}{2k})L(G,W)-\frac{1}{2}{\rm Diag}(We)}{ X}
       \\[5pt]
      {\rm s.t.} & {\rm diag}(X) = e, \;\; X  \in \S^n_+,\;\; X - M \in \cK_{\cP},
 \end{array}
\end{eqnarray}
where $\cK_{\cP} = \{X \in \cS^n \mid X_{ij} = 0, \forall (i,j) \in U; X_{ij} \geq  0, \forall (i,j) \in E \setminus U \} $.

We should mention that we can easily extend our algorithm to handle the following slightly  more general doubly non-negative SDP:
\begin{eqnarray}
  \max \left\{ \inprod{-C}{X} \mid\cA_E X=  b_{ E}, \; X \in \cS_+^n, \; X - M \in \cK_p \right\},  \nonumber
\end{eqnarray}
where $ M \in \S^n$ is a  given  matrix. Thus \eqref{eq-fap-sdpp} can also be
solved by our proposed  algorithm.

\subsubsection{Numerical results}
For large-scale DNN-SDP problems, there exist two  other competitive codes in the literature that are based on  alternating direction   algorithms: a directly extended ADMM  solver (called SDPAD  in \cite{wen2010alternating}) and a two-easy-block-decomposition hybrid proximal extragradient method solver (called 2EBD-HPE in \cite{monteiro2013first} but we will just call it as 2EBD for convenience).

 Here we compare our algorithm ADMM3c with SDPAD (release-beta2, released in December 2012), and 2EBD\footnote{\url{www2.isye.gatech.edu/~cod3/CamiloOrtiz/Software_files/2EBD-HPE_v0.2/2EBD-HPE_v0.2.zip}} (v0.2, released on May 31, 2013)
for solving DNN-SDP. We also include
a convergent alternating direction method with a Gaussian back substitution   proposed in \cite{he2012alternating} (we call the method ADMM3g here and
use the parameter $\alpha=0.999$ in the Gaussian back substitution step\footnote{We avoid taking $\alpha =1$ as it leads to slow convergence for  quite a number of  tested examples.}). We have implemented both  ADMM3c and ADMM3g in {\sc Matlab}. The computational results for all the DNN-SDP problems are obtained on a Linux server (6-core, Intel Xeon X5650 @ 2.67GHz, 32G RAM).

We measure the accuracy of an approximate optimal solution $(X,y_E,S,Z)$ for (\ref{eq-primal-DNNSDP}) and (\ref{eq-dual-DNNSDP}) by using the following relative residual:
\begin{eqnarray}
\eta = \max \{\eta_P,\eta_D,\eta_{\cK},\eta_{\cP},\eta_{\cK^*},\eta_{\cP^*},\eta_{C_1},\eta_{C_2}\},  \label{eq-eta}
\end{eqnarray}
where
\begin{eqnarray}
&\eta_P =\frac{\norm{\cA_E X-b_E}}{1+\norm{b_E}},\;
\eta_D =\frac{\norm{\cA_E^* y_E+S +Z-C}}{1+\norm{C}},\;
\eta_{\cK} = \frac{\norm{\Pi_{\cS_+^n}(-X)}}{1+\norm{X}},\;
\eta_{\cP} = \frac{\norm{X-\Pi_{\cK_\cP}(X)}}{1+\norm{X}},&
\nonumber \\[5pt]
&\eta_{\cK^*} = \frac{\norm{\Pi_{\cS_+^n}(-S)}}{1+\norm{S}},\;
\eta_{\cP^*} = \frac{\norm{Z-\Pi_{\cK_{\cP}^*}(Z)}}{1+\norm{Z}},\;
\eta_{C_1} = \frac{\abs{\inprod{X}{S}}}{1+\norm{X}+\norm{S}},\;
\eta_{C_2} = \frac{\abs{\inprod{X}{Z}}}{1+\norm{X}+\norm{Z}}.&
\label{eq-etasub}
\end{eqnarray}
Additionally, we compute the relative gap by
\begin{eqnarray}
&\eta_g = \frac{\inprod{C}{X} - \inprod{b_E}{y_E} }{1+\abs{\inprod{C}{X}}+\abs{\inprod{b_E}{y_E}}}.&
\label{eq-gap}
\end{eqnarray}
We terminate the solvers ADMM3c, ADMM3g and SDPAD when
$\eta  < 10^{-6}. $
Note that, as mentioned in the introduction, the direct extension of ADMM to the case of a multi-block problem is not necessarily convergent \cite{chen2013direct}. Hence SDPAD, which is essentially an implementation of ADMM3d with $\tau =1.618$ for solving DNN-SDP problems,  does not have  convergence guarantee. For the implementation of 2EBD, we need to explain in more details.

The method 2EBD in \cite{monteiro2013first} is
designed to solve a conic optimization problem of the form
\begin{eqnarray}
  \min \Big\{ \inprod{C}{X} \mid \cA_1(X) - b_1 \in \cC_1, \;\cA_2(X) - b_2 \in \cC_2 \Big\},
  \label{eq-2EBD-primal}
\end{eqnarray}
where $b_1 \in \cW_1$, $b_2 \in \cW_2$, $C\in \cX$ are given data, $\cA_1:\cX\rightarrow \cW_1$, $\cA_2:\cX\rightarrow \cW_2$ are given linear maps, and $\cC_1 \in \cW_1$, $\cC_2 \in \cW_2$ are nonempty closed convex cones.
The dual of (\ref{eq-2EBD-primal}) is given by
\begin{eqnarray}
 \max \Big\{ \inprod{b_1}{w_1} + \inprod{b_2}{w_2} \mid \cA_1^*(w_1) + \cA_2^*(w_2) = C,\;
  w_1\in \cC_1^*, \; w_2\in\cC_2^* \Big\}.
  \label{eq-2EBD-dual}
\end{eqnarray}
Note that the application of 2EBD to DNN-SDP strongly depends on the possibility of splitting its constraints into two-easy blocks such that the projection problems  $\min \big\{ \frac{1}{2}\norm{X-X_0}^2 \mid \cA_i(X) - b_i \in \cC_i \big\}, i = 1, 2$ can be   computed easily
for any given ${X}_0$.  The users need to input the algorithms for computing $\min \big\{ \frac{1}{2}\norm{X-X_0}^2 \mid \cA_i(X) - b_i \in \cC_i \big\}, i = 1, 2$.
For BIQ, $\theta_{+}$ and FAP problems, their constraints can naturally be split into two-easy blocks \cite{monteiro2013first}.
For  general DNN-SDP problems, however it may be difficult to split their constraints naturally into  two-easy blocks although one can always   reformulate them as SDP problems in the standard form. In our numerical experiments, we use this approach to test   QAP and RCP problems as these problems do not appear to have obvious two-easy blocks structures.

For testing 2EDB on DNN-SDP problems, in order to agree with our stopping criterion, given a solution $(X,{w}_1,{w}_2,S,Z)$ for (\ref{eq-2EBD-primal}) and (\ref{eq-2EBD-dual}), where $S, Z$ are the dual variables corresponding to $X \in \cS_+^n$ and $X \in \cK_{\cP}$  respectively, we measure the relative residual as follows:
\begin{eqnarray}
\hat \eta = \max \{\eta_P,\eta_D,\eta_{\cK},\eta_{\cP},\eta_{\cK^*},\eta_{\cP^*},\eta_{C_1},\eta_{C_2}\}  \label{eq-2EBD-eta},
\end{eqnarray}
where $\eta_{\cK}, \eta_{\cP}, \eta_{\cK^*}, \eta_{\cP^*}, \eta_{C_1},
\eta_{C_2}$ are defined as in \eqref{eq-etasub}, and
$\eta_P =\frac{\norm{(\Pi_{\cC_1^*}(b_1-\cA_1 X),\, \Pi_{\cC_2^*}(b_2-\cA_2 X))} }{1+\norm{(b_1, b_2)}}$,
$\eta_D =\frac{\norm{\cA_1^*(w_1) + \cA_2^*(w_2) - C}}{1+\norm{C}}$.
We terminate the solver 2EBD when
$
\hat \eta < 10^{-6}.
$
And we measure the relative gap as
\begin{eqnarray}
\hat \eta_g = \frac{\inprod{C}{X} - (\inprod{b_1}{w_1} + \inprod{b_2}{w_2})}{1+\abs{\inprod{C}{X}}+\abs{\inprod{b_1}{w_1} + \inprod{b_2}{w_2}}}.
\end{eqnarray}

{We should mention in the implementations of all the solvers, ADMM3c, SDPAD, ADMM3g, and
2EDB, the penalty parameter $\sigma$ is dynamically adjusted according to the progress
of the algorithms. In addition, all the algorithms also adopt some kind of restart strategies
to ameliorate slow convergence.
The exact details on the adjustment strategies are too tedious to
be presented here but it suffices to mention that the key idea is to adjust $\sigma$ so as to
balance the progress of primal  feasibilities $(\eta_P,\eta_{\cK},\eta_{\cP})$ and
dual feasibilities $(\eta_D,\eta_{\cK^*},\eta_{\cP^*})$.
In our numerical experiments, we use the same adjustment strategy for both the solvers
ADMM3c and
ADMM3g to solve all the tested problems, i.e., we do not change the strategy
to maximize the performance of different classes of tested problems.
The solver SDPAD also uses a common adjustment strategy, though different from
that of ADMM3c and ADMM3g. But 2EDB uses
different parameter settings for the adjustment strategy for different
classes of tested problems.

Table \ref{table-numsolved} shows the number  of problems that  have been successfully solved to the accuracy of  $10^{-6}$ in $\eta$ or $\hat{\eta}$ by each of  the four solvers ADMM3c, SDPAD,  ADMM3g and 2EBD, with the  maximum number of iterations set at $25000$.
We can see that ADMM3c  solved the most number of instances
to the required accuracy, with SDPAD in the second place, followed by ADMM3g and 2EDB in the third and fourth place, respectively.
It is actually quite remarkable that all the four solvers are able
to solve these large scale SDP problems to such a good accuracy
despite the fact that they are all first order methods.

Table \ref{tableLS_CADM_SDPAD_ADMG_2EBD}  reports detailed numerical results
for ADMM3c, SDPAD, 2EBD and ADMM3g in solving
some very large scale  DNN-SDP problems.
The detailed results for all the $414$ problems  tested can be found
at {\tt http://www.math.nus.edu.sg/$\sim$mattohkc/publist.html/}.
 Note that we did not list the numerical results for the directly extended ADMM with $\tau=1$ here as it  almost always takes  $20\%$ to $50\%$ more time than the one  with $\tau =1.618$, i.e., SDPAD.
From the detailed numerical results,
one can observe that ADMM3c is generally the fastest in terms of the
computing time, especially when the problem size  is large.

\begin{center}
\begin{longtable}{| c | c | c | c | c |}
\caption{Numbers of problems which are solved to the accuracy of  $ 10^{-6}$
in $\eta$ or $\hat{\eta}$.} \label{table-numsolved}
\\
\hline
\multicolumn{1}{|c}{problem set (No.) $\backslash$ solver} & \multicolumn{1}{|c|}{ADMM3c} & \multicolumn{1}{|c|}{SDPAD}   & \multicolumn{1}{|c|}{2EBD} & \multicolumn{1}{|c|}{ADMM3g} \\ \hline
\endhead
\hline
\endfoot
$\theta_{+}$ (58) & 58 & 58 & 56 & 54 \\ \hline FAP ( 7) &  7 &  7 &  7 &  7 \\ \hline QAP (95) & 39 & 30 & 16 & 28 \\ \hline BIQ (134) & 134 & 134 & 134 & 130 \\ \hline RCP (120) & 120 & 114 & 109 & 113 \\ \hline Total (414) & 358 & 343 & 322 & 332 \\ \hline 
\hline
\end{longtable}
\end{center}

Figure \ref{fig-all} shows the performance profiles in terms of computing time for  ADMM3c, SDPAD, 2EBD and ADMM3g, for all the tested problems including  those problems not listed in Table \ref{tableLS_CADM_SDPAD_ADMG_2EBD}. We recall that a point $(x, y)$ is in the performance profiles curve of a method if and only if it can solve $(100y)\%$ of all the tested problems no slower than $x$ times of any other methods. It can be seen that ADMM3c outperforms all the other 3 solvers by a significant margin.

Figure \ref{fig-all-iter} shows the performance profiles in terms of the number of iterations for ADMM3c, SDPAD, 2EBD and ADMM3g, for all the tested problems. We may observe that for
the majority of  the tested problems, ADMM3c takes the least number of iterations.
For the BIQ problems, the solver 2EDB typically takes the least number of iterations.
However, as each iteration of 2EDB requires quite a number of intermediate
calculations to estimate a step-length to achieve good convergence, the
non-trivial overheads incurred often counteract the savings in the number of
iterations. As a result, even though the performance
profile of 2EDB in terms of the number of iterations  dominates that of ADMM3g,
its profile in terms of the computing time
does not behave similarly.

\begin{figure}[t]
\begin{tabular}{cc}
 \begin{minipage}[t]{0.5\textwidth}
 \ifpdf
 \includegraphics[width=0.95\textwidth]{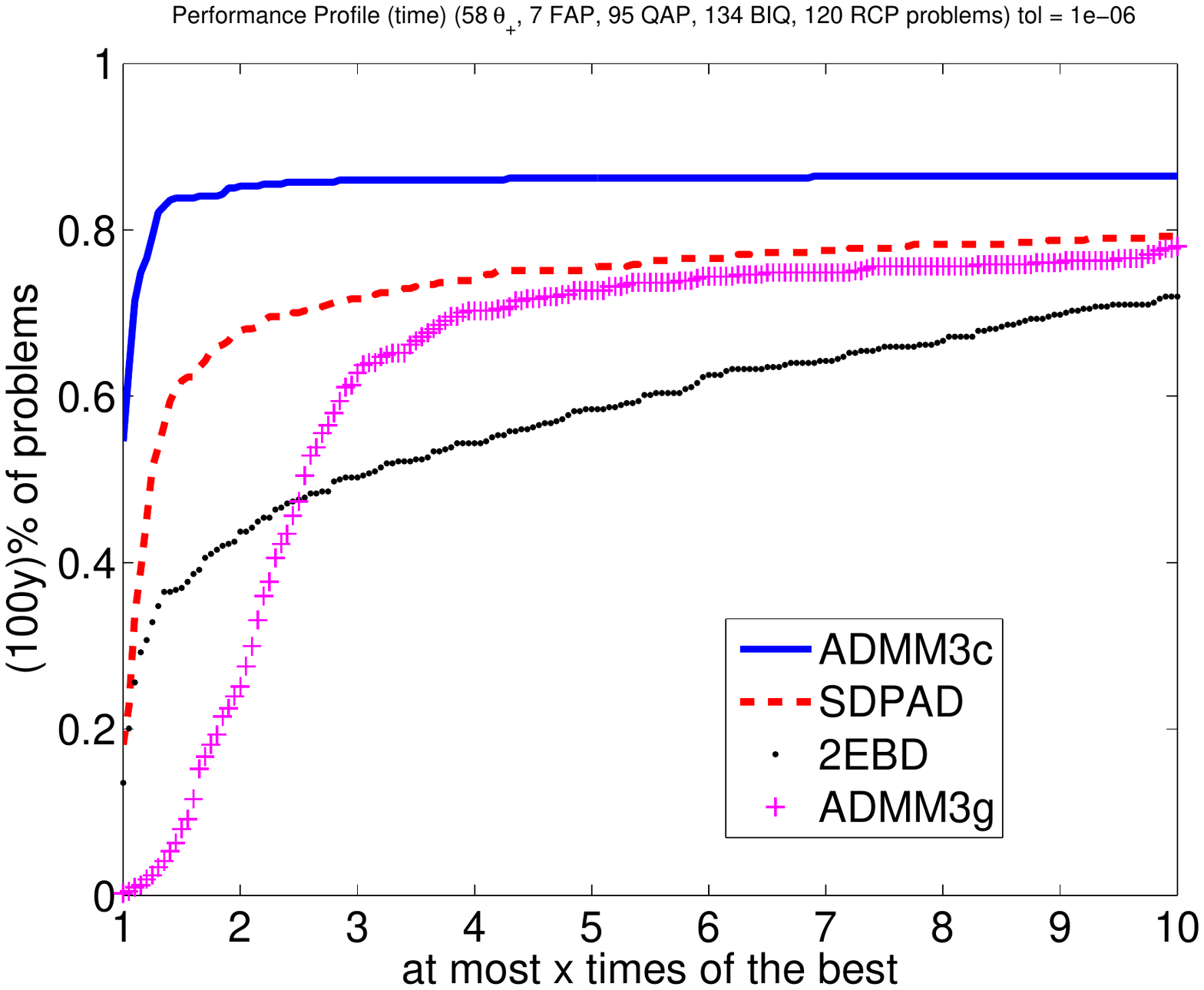}
 \else
 \includegraphics[width=0.95\textwidth]{Performance_Profile_alltime.eps}
 \fi
 \caption{\footnotesize{Performance profiles (time) of ADMM3c, \;\;\;
SDPAD,  ADMM3g and 2EBD.}}
 \label{fig-all}
 \end{minipage}

 \begin{minipage}[t]{0.5\textwidth}
 \ifpdf
 \includegraphics[width=0.95\textwidth]{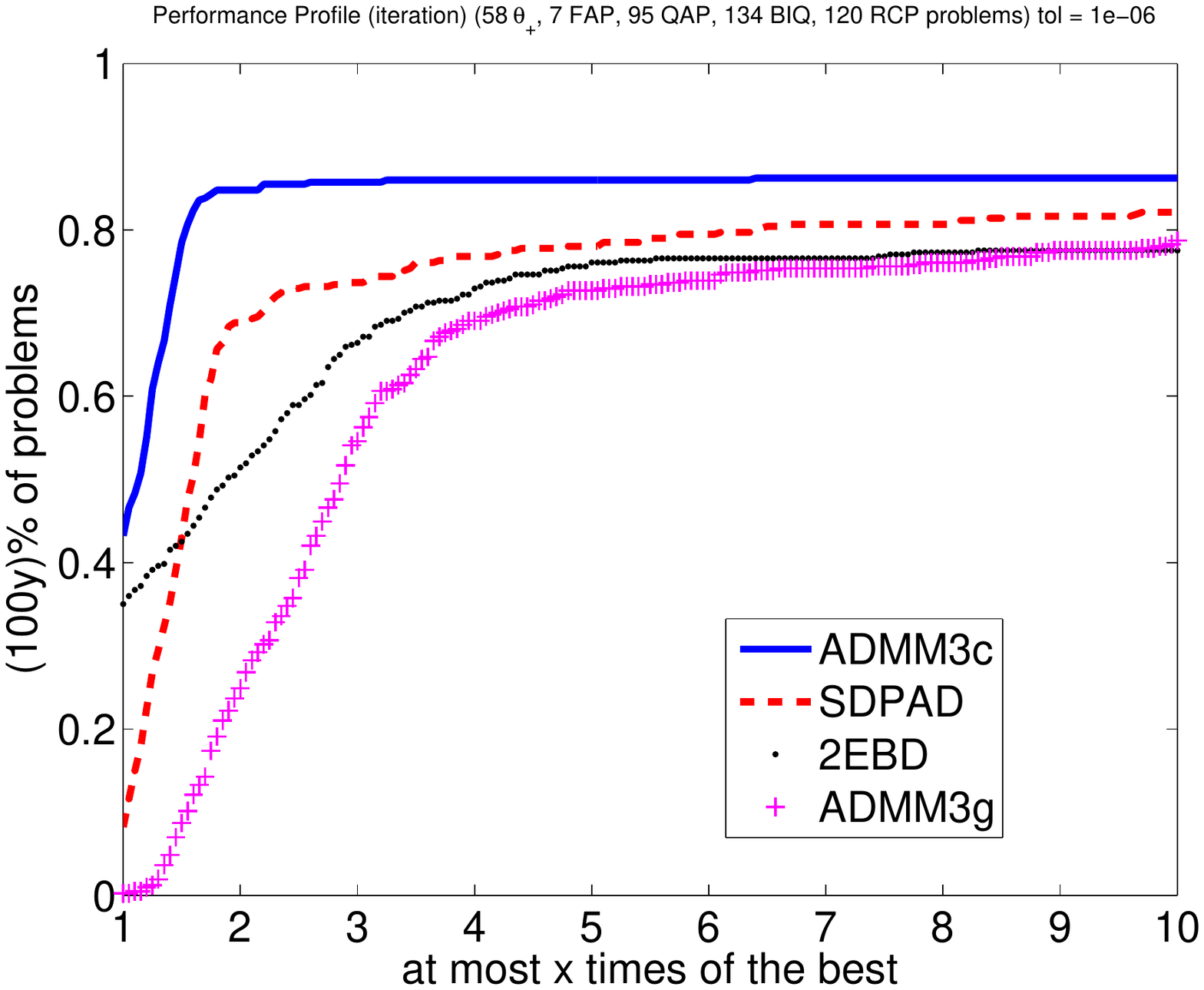}
 \else
  \includegraphics[width=0.95\textwidth]{Performance_Profile_alliteration.eps}
  \fi
 \caption{\footnotesize{Performance profiles (iteration) of ADMM3c, SDPAD, ADMM3g and 2EBD.}}
 \label{fig-all-iter}
\end{minipage}
\end{tabular}
\end{figure}

\subsection{Numerical results for SDP  with  many inequality constraints}\label{subsection:SDPInequality}

In this subsection, we will consider
  {\bf (SDP)} with many  inequality constraints $\cA_I x_I \geq b_I$.
The dual of {\bf (SDP)} takes the form of
  \begin{eqnarray}
  \min \left\{ \delta _{{\cS}^n_+} (S) +( \delta _{\Re_+^{m_I}} (y_I)    -\inprod{b_I}{y_I})+\delta_{\cK^*_p}(Z) - \inprod{b_E}{y_E}   \mid S+\cA_I^* y_I +Z+ \cA^*_E y_E = C
     \right\}.
  \label{eq-dualSDP}
\end{eqnarray}
 Let $\cD:\cS^n \to \cS^n$ be a   given nonsingular linear operator    and $\cD^*$ be its adjoint. By introducing an extra variable, we can reformulate (\ref{eq-dualSDP}) into the form of (\ref{eq-dual01-AoA}), for which our proposed Conic-sPADMM3c can be used.

For the SDP problems described in (\ref{eq-BIQ-SDP}) arising from relaxing the BIQ problems, in order to get tighter bounds, we may add  in some valid inequalities to  get the following problems:
\begin{eqnarray}\begin{array}{ll}
   \min & \frac{1}{2} \inprod{Q}{Y} + \inprod{c}{x} \\[5pt]
  {\rm s.t.} &  {\rm diag}(Y) - x = 0, \quad \alpha = 1, \quad
   X = \left[ \begin{array}{cc} Y & x \\ x^T & \alpha  \end{array}\right]\in \S^n_+,\quad
  X \in \cK_{\cP},
\\[8pt]
  & -Y_{ij} + x_i \geq 0,\;\;
 -Y_{ij} + x_j \geq 0,\;\;
 Y_{ij} - x_i -x_j \geq -1,\; \forall \; i < j, \; j=2,\dots,n-1,
  \end{array}
  \label{eq-BIQineq-SDP}
\end{eqnarray}
where $\cK_{\cP} = \{X \in \cS^n \mid X \geq 0 \}$. For convenience, we call the problem in (\ref{eq-BIQineq-SDP}) an extended BIQ problem.
Note that the last set of inequalities in \eqref{eq-BIQineq-SDP}
are obtained from the valid inequalities $x_i(1-x_j)\geq 0$, $x_j(1-x_i)\geq 0$, $(1-x_i)(1-x_j)\geq 0$
when $x_i,x_j$ are binary variables.

Note that one may also apply the directly  extended ADMM (see (\ref{admmlb}))  with 4 blocks by adding a proximal term for the $y_I$-part (similar to the discussion in Section \ref{subsection:ConicD-4block}). We call this method {\sc sPadmm}4d. Of course, we are mindful that {\sc sPadmm}4d has no convergence guarantee.
In this subsection, we compare the algorithms {\sc sPadmm}3c, {\sc sPadmm}4d and a convergent linearized alternating direction method with a Gaussian back substitution  proposed in \cite{HYuan2013} (we
 call the method {\sc Ladmm}4g here and use the parameter $\alpha = 0.999$ in the Gaussian back substitution step)   for the extended BIQ problems (\ref{eq-BIQineq-SDP}). We have implemented {\sc sPadmm}3c, {\sc sPadmm}4d and {\sc Ladmm}4g in {\sc Matlab}.
For {\sc sPadmm}4d, we set the step-length $\tau = 1.618$.
For the purpose of comparison, we also test the directly extended
sPADMM with unit step-length (i.e., $\tau=1$), which is called {\sc sPadmm}4d(1). The computational results for all the extended BIQ problems are obtained on the same   Linux server as before.

We note here that neither SDPAD nor 2EBD can be directly applied to solve  the
problems \eqref{eq-BIQineq-SDP}. One may of course  try to first convert the inequality constraints into linear equalities by introducing slack variables and then apply both SDPAD and 2EBD to the reformulated problems. However, such an approach is inefficient as the
linear system of equations which needs to be solved at each iteration is
very large but not-so-sparse.
Not surprisingly, this approach is  very slow for the extended BIQ problems according to our  numerical experience.

We measure the accuracy of an approximate optimal solution $(X,y_E,y_I,S,Z)$ for {\bf (SDP)} and its dual (\ref{eq-dualSDP}) by using the following relative residual:
\begin{eqnarray}
\eta = \max \{\eta_P,\eta_D,\eta_{\cK},\eta_{\cP},\eta_{\cK^*},\eta_{\cP^*},\eta_{C_1},\eta_{C_2},\eta_{I},\eta_{I^*}\},  \label{eq-eta-4block}
\end{eqnarray}
where $\eta_{\cK}$,
$\eta_{\cP}$,
$\eta_{\cK^*}$,
$\eta_{\cP^*}$,
$\eta_{C_1} $,
$\eta_{C_2} $, are defined as in \eqref{eq-etasub}, and
\begin{eqnarray*}
& \eta_P =\frac{\norm{\cA_E X-b_E}}{1+\norm{b_E}},\;
\eta_D =\frac{\norm{\cA_E^* y_E+\cA_I^* y_I+S +Z-C}}{1+\norm{C}},\;
\eta_I =\frac{\norm{{\rm max}(0,b_I-\cA_I X)}}{1+\norm{b_I}},\;
\eta_{I^*} =\frac{\norm{{\rm max}(0,-y_I)}}{1+\norm{y_I}}. &
\end{eqnarray*}
Additionally, we compute the relative gap by
\begin{eqnarray}
&\eta_g = \frac{\inprod{C}{X} - (\inprod{b_E}{y_E}+\inprod{b_I}{y_I}) }{1+\abs{\inprod{C}{X}}+\abs{\inprod{b_E}{y_E}+\inprod{b_I}{y_I}}}.&
\label{eq-gap-4block}
\end{eqnarray}
We terminate {\sc sPadmm}3c, {\sc sPadmm}4d, {\sc sPadmm}4d(1) and {\sc Ladmm}4g when
$
\eta  < 10^{-5}
$ or when they reach the maximum number of $50000$ iterations.

In Table \ref{tableLS_CADM_ADMG_ADM_ADM0}, we   report some detailed numerical results for the solvers {\sc sPadmm}3c, {\sc Ladmm}4g, {\sc sPadmm}4d and {\sc sPadmm}4d(1) in solving a collection of $134$ extended BIQ problems.

Figure \ref{fig-biqineq} shows the performance profiles  in terms of computing time for
{\sc sPadmm}3c, {\sc Ladmm}4g, {\sc sPadmm}4d and {\sc sPadmm}4d(1) in solving 134 extended BIQ problems.
One can observe  that {\sc Ladmm}4g is much slower than the other three solvers.
The solver {\sc sPadmm}3c is clearly  more  efficient than the directly extended
sPADMM with unit step-length, i.e., {\sc sPadmm}4d(1),  and it is even faster than {\sc sPadmm}4d with $\tau=1.618$,
though only marginally.

Figure \ref{fig-biqineq-iter} shows the performance profiles in terms of the  number of iterations for {\sc sPadmm}3c, {\sc Ladmm}4g, {\sc sPadmm}4d and {\sc sPadmm}4d(1).
Observe that for the majority of the test problems, {\sc sPadmm}3c takes less iterations than {\sc sPadmm}4d with step-length $\tau = 1.618$. However, for some test problems, due to the
overheads incurred in handling the additional matrix variable introduced to reformulate
\eqref{eq-dualSDP} into the form \eqref{eq-dual01-AoA}, {\sc sPadmm}3c may
take slightly more time than {\sc sPadmm}4d,
even though  the latter may take slightly more iterations.

\begin{figure}[t]
\begin{tabular}{cc}
 \begin{minipage}[t]{0.5\textwidth}
 \ifpdf
 \includegraphics[width=0.95\textwidth]{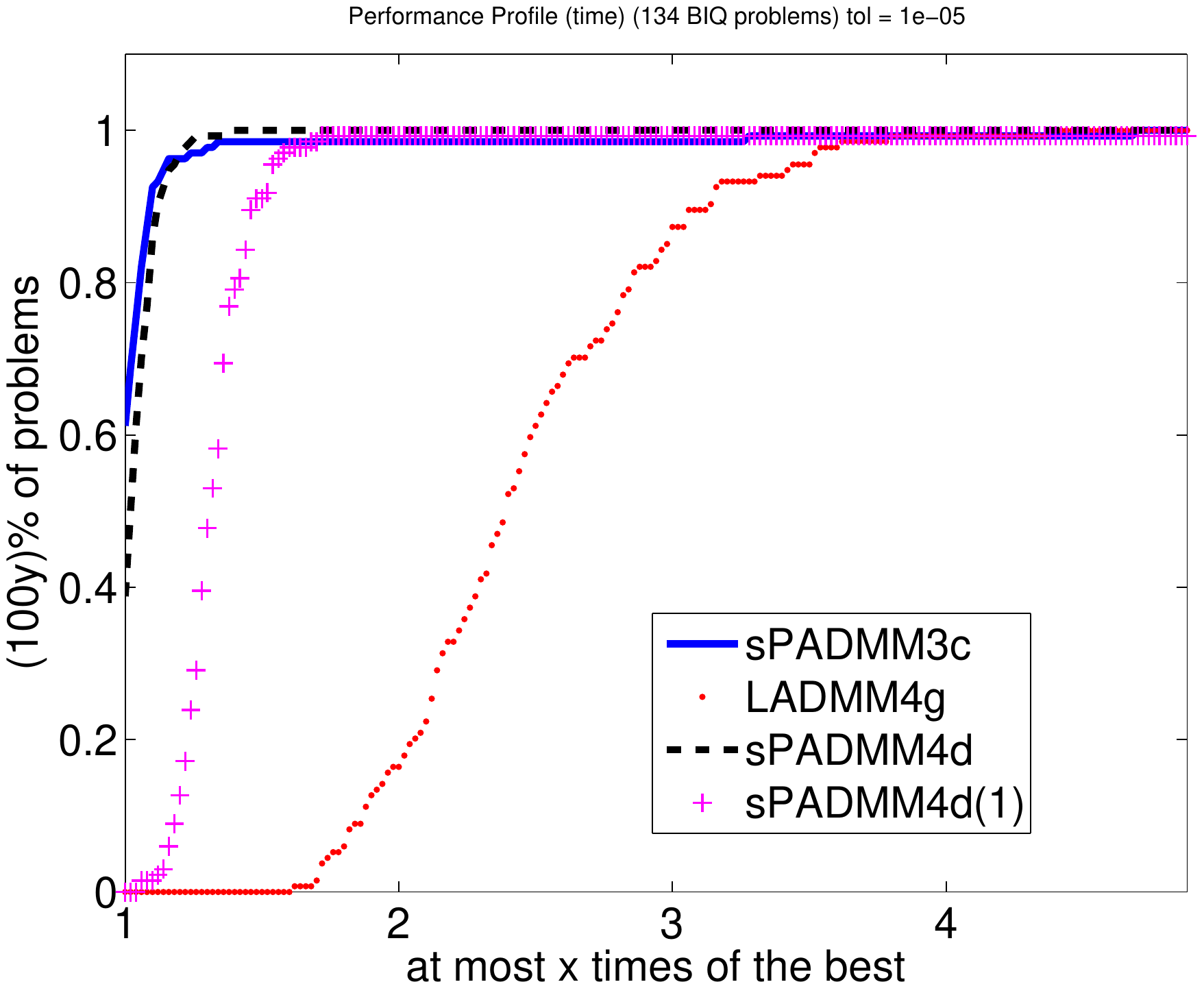}
 \else
 \includegraphics[width=0.95\textwidth]{Performance_Profile_time.eps}
 \fi
 \caption{\footnotesize{Performance profiles (time) of {\sc sPadmm}3c, {\sc Ladmm}4g, {\sc sPadmm}4d and {\sc sPadmm}4d(1).}}
 \label{fig-biqineq}
 \end{minipage}

 \begin{minipage}[t]{0.5\textwidth}
 \ifpdf
 \includegraphics[width=0.95\textwidth]{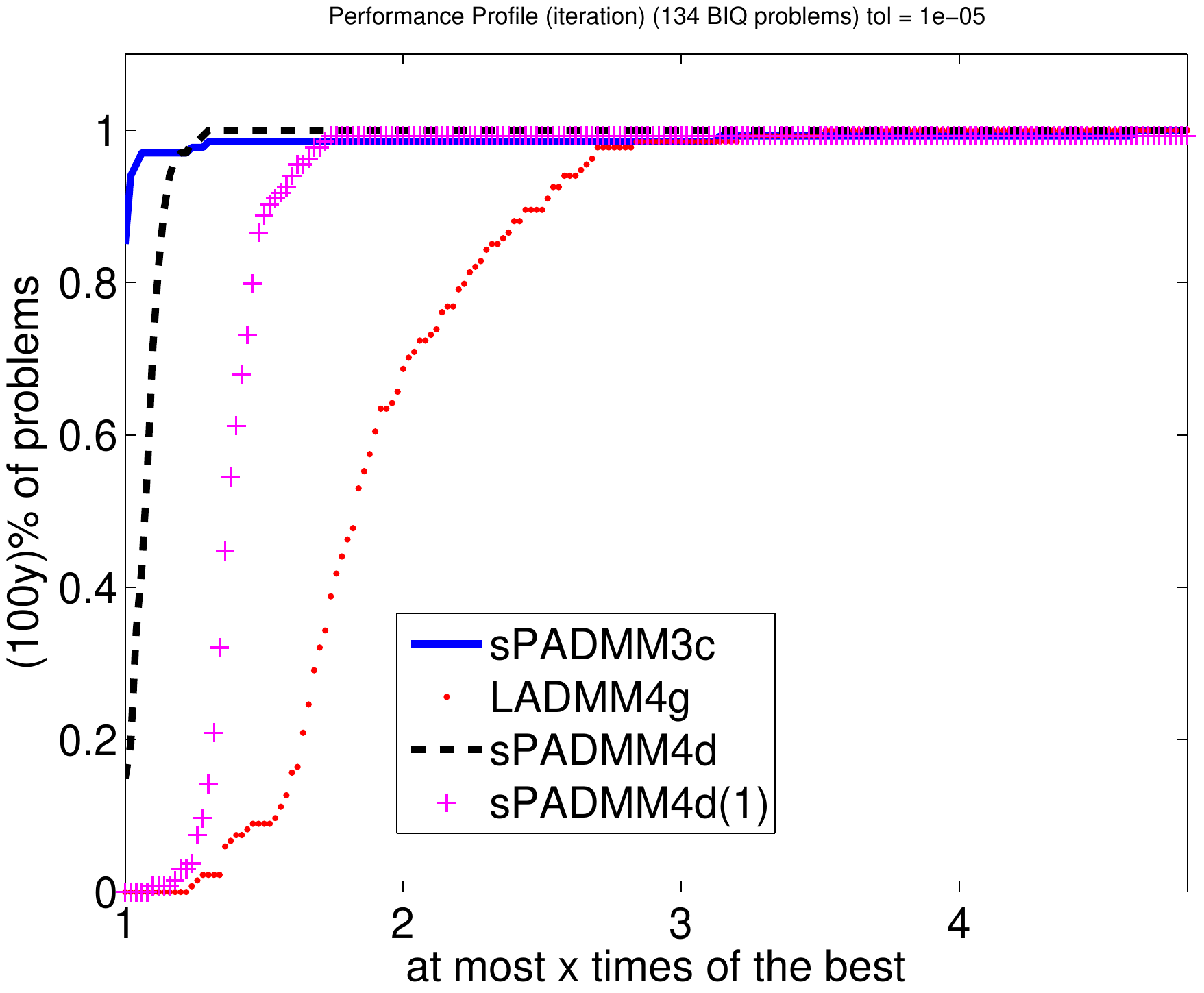}
 \else
  \includegraphics[width=0.95\textwidth]{Performance_Profile_iteration.eps}
  \fi
 \caption{\footnotesize{Performance profiles (iteration) of {\sc sPadmm}3c, {\sc Ladmm}4g, {\sc sPadmm}4d and {\sc sPadmm}4d(1).}}
 \label{fig-biqineq-iter}
\end{minipage}
\end{tabular}
\end{figure}

\section{Conclusions}\label{section:conclusions}

In this paper, we have proposed a 3-block semi-proximal ADMM  that is  both convergent and efficient for  finding a solution of  medium accuracy to  conic programming problems with four types  of constraints.  By conducting numerical experiments on a large number of
   doubly non-negative SDP problems  with equality and/or inequality constraints, we have presented convincing numerical results showing that for the
vast majority of problems tested, our proposed   (semi-proximal) ADMM is   at least  $20\%$   faster than the directly extended (semi-proximal) ADMM with unit step-length.  At least for the class of
conic programming {(\bf P)} problems, we can safely say that we have  resolved the dilemma that an ADMM    is either efficient in practice  but without convergent guarantee in theory or  the contrary.  This   opens up the possibility  of designing  a convergent and yet
practically
efficient  ADMM   with an intelligent  BCD cycle
rather than  the usual non-convergent   Gauss-Seidel BCD cycle for solving multi-block convex optimization problem (\ref{eq-ADMMq}).
 We leave this as one of our future research topics. In fact, our  primary motivation of  introducing this convergent 3-block semi-proximal ADMM  is   to quickly  produce an initial point
for    conic programming   ({\bf P}) so as to  warm-start methods which have fast local convergence properties. For SDP problems in  standard form,
this has already been done by
Zhao, Sun and Toh  in \cite{ZSToh2010} by first using the classic 2-block ADMM (it was called the boundary point method in \cite{povh2006boundary} at that time) to generate a starting point and then switching it to the fast convergent Newton-CG augmented Lagrangian method. The resulting software SDPNAL  has  been successfully employed by Nie and Wang \cite{NWang2012, NWang2013} to  solve very large scale  SDP problems in standard form   arising from polynomial optimization and  rank-1 tensor approximation problems. Naturally, our next target is to extend this approach to multi-block convex optimization problems beyond
conic programming  ({\bf P}) in standard form. We will report our corresponding findings in subsequent works.

\section*{Acknowledgements.}
The authors would like to thank Dr Caihua Chen at Nanjing University and Professor Deren Han at Nanjing Normal University for  bringing our attention to several  papers  on the convergence properties of multi-block ADMM and its various modifications and thank Mr Xudong Li at National University of Singapore for numerous discussions on many   topics covered in this paper. Thanks also go to the three referees for their comments and suggestions that help us  improve the
paper.
\bibliographystyle{siam}
\bibliography{SDPbib}


\begin{landscape}
\begin{scriptsize}
\begin{center}

\begin{scriptsize}

\begin{longtable}{| ccc | c | c | c | c |}
\caption{The performance of ADMM3c, SDPAD, ADMM3g, 2EBD on $\theta_+$, FAP, QAP, BIQ and RCP problems (accuracy $= 10^{-6}$). In the table, ``3c" and ``3g" stand
for ADMM3c and ADMM3g, respectively. The computation time is in the
format of ``hours:minutes:seconds".}
\label{tableLS_CADM_SDPAD_ADMG_2EBD} \\
\hline
\multicolumn{3}{|c}{} & \multicolumn{1}{|c|}{iteration} & \multicolumn{1}{|c|}{$\eta|\,\eta|\,\eta|\,\hat{\eta}$} & \multicolumn{1}{|c|}{$\eta_g$}  & \multicolumn{1}{|c|}{time} \\ \hline
problem & $m_E;m_I$ & $n_{s};$ &3c$|$SDPAD$|$3g$|$2EBD&3c$|$SDPAD$|$3g$|$2EBD&3c$|$SDPAD$|$3g$|$2EBD&3c$|$SDPAD$|$3g$|$2EBD \\ \hline
\endhead

 theta10 &  12470;0 &   500; &  354 $|$  351 $|$  636 $|$  490&  8.5-7  $|$  9.9-7  $|$  9.9-7  $|$  9.9-7  &  {\green{ -2.5-6}}  $|$  -7.8-7  $|$  7.4-7  $|$  {\green{ 1.5-6}}  & 45 $|$ 42 $|$ 1:34 $|$ 1:11  \\ \hline 
 theta102 &  37467;0 &   500; &  157 $|$  130 $|$  232 $|$  355&  9.5-7  $|$  9.0-7  $|$  9.7-7  $|$  9.9-7  &  -6.0-7  $|$  {\green{ 2.1-6}}  $|$  9.0-7  $|$  {\green{ 2.2-6}}  & 23 $|$ 23 $|$ 36 $|$ 54  \\ \hline 
 theta103 &  62516;0 &   500; &  144 $|$  108 $|$  199 $|$  323&  9.2-7  $|$  9.8-7  $|$  9.9-7  $|$  9.8-7  &  -3.0-8  $|$  5.9-7  $|$  1.8-7  $|$  {\green{ 2.6-6}}  & 22 $|$ 21 $|$ 32 $|$ 49  \\ \hline 
 theta104 &  87245;0 &   500; &  169 $|$  123 $|$  226 $|$  338&  9.3-7  $|$  9.8-7  $|$  9.0-7  $|$  9.9-7  &  -9.2-8  $|$  {\green{ 1.1-6}}  $|$  4.2-7  $|$  {\green{ 3.2-6}}  & 24 $|$ 20 $|$ 36 $|$ 51  \\ \hline 
 theta12 &  17979;0 &   600; &  362 $|$  366 $|$  648 $|$  494&  9.0-7  $|$  8.8-7  $|$  9.6-7  $|$  9.2-7  &  {\green{ -2.2-6}}  $|$  {\green{ 1.0-6}}  $|$  {\green{ 1.2-6}}  $|$  {\green{ 1.3-6}}  & 1:14 $|$ 1:11 $|$ 2:28 $|$ 1:52  \\ \hline 
 theta123 &  90020;0 &   600; &  156 $|$  107 $|$  197 $|$  345&  9.3-7  $|$  9.9-7  $|$  9.5-7  $|$  9.9-7  &  -6.0-8  $|$  5.1-7  $|$  1.4-7  $|$  {\green{ 2.6-6}}  & 35 $|$ 33 $|$ 49 $|$ 1:26  \\ \hline 

 c-fat200-1 &  18367;0 &   200; &  233 $|$  444 $|$  472 $|$  330&  9.8-7  $|$  9.9-7  $|$  9.5-7  $|$  9.9-7  &  -6.9-7  $|$  {\green{ -1.2-6}}  $|$  {\green{ -1.1-6}}  $|$  {\green{ 2.1-6}}  & 03 $|$ 06 $|$ 06 $|$ 04  \\ \hline 
 hamming-8-4 &  11777;0 &   256; &  124 $|$  104 $|$  179 $|$  214&  4.7-7  $|$  9.6-7  $|$  8.5-7  $|$  8.9-7  &    {\blue{ -5.3-6}}  $|$  {\green{ 2.1-6}}  $|$  {\green{ -1.4-6}}  $|$    {\red{ 1.0-5}}  & 02 $|$ 03 $|$ 04 $|$ 04  \\
\hline 
 
\end{longtable}

\newpage
\begin{longtable}{| ccc | c | c | c | c |}
\caption{The performance of {\sc sPadmm}3c, {\sc Ladmm}4g, {\sc sPadmm}4d,  {\sc sPadmm}4d(1) on extended BIQ problems (accuracy $= 10^{-5}$). In the table, we have omitted the
command string ``ADMM" in the names of the solvers. The computation time is in the
format of ``hours:minutes:seconds".}
\label{tableLS_CADM_ADMG_ADM_ADM0} \\
\hline
\multicolumn{3}{|c}{} & \multicolumn{1}{|c|}{iteration} & \multicolumn{1}{|c|}{$\eta$} & \multicolumn{1}{|c|}{$\eta_g$}  & \multicolumn{1}{|c|}{time} \\ \hline
problem & $m_E;m_I$ & $n_{s};$ &P3c$|$L4g$|$P4d$|$P4d(1)&P3c$|$L4g$|$P4d$|$P4d(1)&P3c$|$L4g$|$P4d$|$P4d(1)&P3c$|$L4g$|$P4d$|$P4d(1) \\ \hline

\endhead

 be200.3.1 &  201;59700 &   201; &  8943 $|$  14140 $|$  9729 $|$  12236&  9.9-6  $|$  9.9-6  $|$  9.9-6  $|$  9.9-6  &  3.2-7  $|$  -7.1-6  $|$  -1.8-6  $|$  -1.4-7  & 7:34 $|$ 12:45 $|$ 7:22 $|$ 9:04  \\ \hline 
 be200.3.2 &  201;59700 &   201; &  7665 $|$  14170 $|$  8854 $|$  11243&  9.9-6  $|$  9.9-6  $|$  9.9-6  $|$  9.9-6  &  7.8-6  $|$  -7.1-6  $|$  2.4-6  $|$  7.7-6  & 5:27 $|$ 15:12 $|$ 5:54 $|$ 7:56  \\ \hline 
 be200.3.3 &  201;59700 &   201; &  11187 $|$  22332 $|$  11429 $|$  14882&  9.9-6  $|$  9.9-6  $|$  9.9-6  $|$  9.9-6  &  8.9-7  $|$  -7.0-6  $|$  2.4-7  $|$  -7.2-7  & 8:58 $|$ 20:21 $|$ 9:29 $|$ 10:50  \\ \hline 
 be200.3.4 &  201;59700 &   201; &  9513 $|$  21265 $|$  9628 $|$  12833&  9.9-6  $|$  9.9-6  $|$  9.9-6  $|$  9.9-6  &  1.1-6  $|$  -7.1-6  $|$  -1.7-6  $|$  -1.3-6  & 7:11 $|$ 21:43 $|$ 6:58 $|$ 9:59  \\ \hline 
 be200.3.5 &  201;59700 &   201; &  9373 $|$  15368 $|$  9370 $|$  12526&  9.9-6  $|$  9.9-6  $|$  9.9-6  $|$  9.9-6  &  -4.3-7  $|$  {\green{ -1.2-5}}  $|$  -2.9-6  $|$  -9.7-7  & 8:43 $|$ 14:15 $|$ 8:21 $|$ 10:14  \\ \hline 
 be200.3.6 &  201;59700 &   201; &  8777 $|$  13901 $|$  10180 $|$  12967&  9.9-6  $|$  9.8-6  $|$  9.9-6  $|$  9.9-6  &  -5.0-6  $|$  {\green{ -1.1-5}}  $|$  -9.6-6  $|$  -7.3-6  & 7:07 $|$ 14:23 $|$ 7:53 $|$ 9:22  \\ \hline 
 be200.3.7 &  201;59700 &   201; &  12965 $|$  24430 $|$  13904 $|$  17684&  9.9-6  $|$  9.9-6  $|$  9.9-6  $|$  9.9-6  &  -2.6-6  $|$  -8.2-6  $|$  -3.5-6  $|$  -3.7-6  & 10:52 $|$ 22:31 $|$ 12:29 $|$ 14:15  \\ \hline 
 be200.3.8 &  201;59700 &   201; &  11102 $|$  18030 $|$  10906 $|$  14622&  9.9-6  $|$  9.9-6  $|$  9.9-6  $|$  9.9-6  &  3.7-7  $|$  -9.7-6  $|$  -8.5-7  $|$  -1.6-6  & 8:24 $|$ 19:35 $|$ 7:46 $|$ 11:03  \\
\hline

\end{longtable}
\end{scriptsize}

\end{center}
\end{scriptsize}
\end{landscape}

\end{document}